\documentclass[a4,11pt]{article}
\usepackage{subcaption}
\usepackage{graphicx}
\usepackage{epstopdf}
\usepackage{amsmath,amsthm,amsfonts,amssymb,epsfig,graphicx} 
\usepackage{hyperref}
\usepackage{comment}
\usepackage{color,xcolor}
\usepackage[linesnumbered,ruled,vlined]{algorithm2e}
\usepackage{subcaption}

\newcommand{\RR}{{\mathbb{R}}}

\newcommand{\1}{\mbox{\bf 1}}

\newcommand{\bu}{\mbox{\bf u}}
\newcommand{\bv}{\mbox{\bf v}}
\newcommand{\bl}{\mbox{\bf l}}
\newcommand{\br}{\mbox{\bf r}}

\newcommand{\R}{\mathbb{R}}

\newcommand{\rank}{\mathrm{rank}}

\newcommand{\minuh}{\mathrm{min}_{\hat u}}
\newcommand{\minvh}{\mathrm{min}_{\hat v}}
\newcommand{\maxuh}{\mathrm{max}_{\hat u}}
\newcommand{\maxvh}{\mathrm{max}_{\hat v}}

\DeclareMathOperator*{\argmin}{arg\,min}

\newtheorem{theorem}{Theorem}[section]

\newtheorem{prop}[theorem]{Proposition}

\newtheorem{lemma}[theorem]{Lemma}
\newtheorem{corollary}[theorem]{Corollary}
\newtheorem{remark}[theorem]{Remark}
\newtheorem{example}[theorem]{Example}

\definecolor{brilliantrose}{rgb}{1.0, 0.33, 0.64}
\definecolor{myviolet}{rgb}{0.21, 0.0, 0.85}
\definecolor{amethyst}{rgb}{0.6, 0.4, 0.8}
\definecolor{carrotorange}{rgb}{0.93, 0.57, 0.13}

\begin{document}

\title{On rank-2 Nonnegative Matrix Factorizations\\ 
and their variants.}

\author{Etna Lindy, Vanni Noferini and Paul Van Dooren}

\maketitle


\begin{abstract}
We consider the problem of finding the best nonnegative rank-2 approximation of an arbitrary nonnegative matrix. We first revisit the theory, including an explicit parametrization of all possible nonnegative factorizations of a nonnegative matrix of rank $2$. Based on this result, we construct a cheaply computable (albeit suboptimal) nonnegative rank-2 approximation for an arbitrary nonnegative matrix input. This can then be used as a starting point for the Alternating Nonnegative Least Squares method to find a nearest approximate nonnegative rank-2 factorization of the input; heuristically, our newly proposed initial value results in both improved computational complexity and enhanced output quality. We provide extensive numerical experiments to support these claims. Motivated by graph-theoretical applications, we also study some variants of the problem, including matrices with symmetry constraints.
\end{abstract}

\section{Introduction} 
Recall \cite{Gillis} that a real matrix $N \in \R^{m \times n}$ is called nonnegative (and we write $N \geq 0$) if all its entries are nonnegative real numbers; in this case, the \emph{nonnegative rank} of $N$ is the least integer $r$ such that $N=L R^\top$ where $L  \in \R^{m \times r}$, $R \in \R^{n \times r}$, and $L \geq 0$, $R \geq 0$. In the special case where a nonnegative matrix $N$ does not have any zero entry, we write $N>0$ and say that $N$ is positive. Similarly, when a real symmetric matrix $M=M^\top$ is positive semidefinite (resp. positive definite), we write
$M\succeq 0$ (resp. $M\succ 0$).

The task of computing a \emph{Nonnegative Matrix Factorization} (NMF) of a given nonnegative rank consists, given $N$, of computing $L$ and $R$ such that $N=L R^\top$. A related problem is that of computing a rank-$r$ approximate NMF, sometimes also called a \emph{nonnegative approximation}: Given an input $N$ and a nonnegative rank $r$, and having fixed a matrix norm $\| \cdot \|$, we seek nonnegative matrices $L,R$, both with $r$ columns, and such that $\| N - LR^\top \|$ is as small as possible. These exact or approximate NMFs appear naturally in several problems where classes of objects have to be discovered in sets of measurements. Examples range from clustering of hyperspectral images \cite{GilKP14} and mass spectrometry imaging data \cite{MelWM24} to blind separation of signal sources \cite{MouBI05}, or to role modeling of data sets represented as graphs \cite{Barbarino,Marchand}; see also \cite{Gillis} for a nice survey. In this paper we focus on the case of nonnegative rank two. Rank-$2$ nonnegative factorization problems occur, for example, in applied graph theory; in particular, they arise in the context of role modeling when only two ``roles" are being separated \cite{Barbarino,Marchand,RW}, such as in the Core-Periphery problem \cite{CP}. In some related factorization problems of symmetric matrices, one also imposes symmetry \cite{KG12} of the obtained factorization.

Beyond applications, there are also mathematical reasons to pay special attention to rank-$2$ nonnegative matrices. If we denote by $\rank_+$ and $\rank$, respectively, the nonnegative rank and the usual linear-algebraic rank over  the field $\R$, then it is immediate from the definitions that $\rank~N \leq \rank_+~N$ for every nonnegative matrix $N$. However, it has also been shown \cite{MouBI05,Tho74} that equality always holds when $\rank~N \leq 2$. Conversely, if $\rank~N \geq 3$, then there exist examples where the inequality is strict. As a consequence, the case where $N \geq 0$ has linear-algebraic rank at most $2$ stands out, and is more accessible to mathematical analysis, because we can conclude that a factorization $N=L R^\top$ holds,
where both $L$ and $R$ have only two columns, and are nonnegative. The borderline case of rank $2$ is particularly interesting from a complexity theory viewpoint. Indeed, while computing an approximate NMF is NP-hard for $\rank_+(L R^\top) \geq 3$ \cite{Vavasis} and can be done in\footnote{Here and throughout, we use Knuth's $\Theta$ notation to denote asymptotic complexity.} $\Theta(mn)$ flops for $\rank_+(L R^\top)=1$, it is to our knowledge an open problem to determine whether or not the case $\rank_+(L R^\top) = 2$ is NP-hard. Among the above mentioned applications of rank-$2$ approximate NMFs, random graphs with a stochastic block model \cite{Barbarino,Marchand} whose expected adjacency matrix have rank $2$ also motivate us to consider ``three way factorizations" \cite{DacheVG24, Gillis} of rank-$2$ nonnegative matrices. Three way factorizations could be seen as nonnegative generalizations of the $LDL^\top$ decomposition \cite{GolubV,asna} of symmetric matrices.

When studying approximate NMFs, we will focus on the Frobenius norm $\| N \|_F = \sqrt{\sum_{i,j} N_{ij}^2}
$. We will use uppercase for matrices
and reserve the letters $N$, $P$, $L$ and $R$ for, respectively, nonnegative matrices, 
permutation matrices, and left and right matrix factors; other uppercase letters may represent either matrices not in these classes or (less frequently) matrices in these classes when necessary to avoid cumbersome notation.  We use lowercase and bold font for generic vectors, lowercase for generic scalars, lowercase Greek letters for angles, and (bold font) uppercase Greek letters for vectors whose entries are angles. Throughout the paper, we often study closed or open subintervals of $\R$, which we denote respectively by $[a,b]$ or $]a,b[$ where $a \leq b$ are the endpoints of the interval.

The structure of the paper is as follows: In Section \ref{sec:theory} we revisit the theory of rank $2$ NMF, and link it to the Singular Value Decomposition (SVD) of the same input matrix. Then, Theorem \ref{th:main} develops a full parametrization of all rank-$2$ NMFs, which is generalized to ``three way" rank-$2$  NMFs in Theorem \ref{th:threeway}. In Section \ref{sec:algorithm}, we develop a cheap, albeit non-optimal, algorithm for nonnegative rank-$2$ approximation based on a geometric interpretation of Theorem \ref{th:main}. In Section \ref{sec:ANLS}, the rough approximation of Section \ref{sec:algorithm} is then used as a starting point for the Alternating Nonnegative Least Squares (ANLS) method \cite{Gillis}, which is the current best general-purpose algorithm for approximate NMFs. This leads to a significant improvement, measured both in terms of output quality and computational cost, with respect to more standard implementations of ANLS; we illustrate this fact in Section \ref{sec:numexp} by providing numerical experiments. Finally, we draw some conclusions in Section \ref{sec:conclusions}.

\section{Exact rank-$2$ NMF}\label{sec:theory}

In this section, we provide a theoretical analysis of rank-2 NMF. Some of the ideas presented below can be traced back in the literature, e.g., \cite{BG08,KG12,Paatero,VanGZD16}, but our treatment goes beyond what we could find in these references. In particular, much of this section can be considered as an extension to an arbitrary nonnegative matrix of the results in \cite{KG12}, that focuses on the case where $N=N^\top \succeq 0$ is not only nonnegative but also completely positive, i.e., it admits a symmetric NMF $N=LL^\top$. (For a general symmetric nonnegative matrix, complete positivity is a strictly stronger requirement than positive semidefiniteness \cite{Diananda,GW}; but for rank-$2$ matrices, equivalence holds \cite{BermanP}.)

We are given an $m\times n$ nonnegative matrix $N$ of rank $r$. We assume further that the matrix $N$ does not contain any zero row or column. This is no loss of generality for the computation of a NMF; indeed, suppose that $N$ has a zero column (the case of a zero row can be treated analogously). Then, for some permutation matrix $P$, $N=\begin{bmatrix}
    \widetilde N & {\bf 0}
\end{bmatrix} P$. Assume that we can compute a NMF $\widetilde N = L\widetilde R^\top$, $L,\widetilde R \geq 0$. Then, $N = L \left(\begin{bmatrix}
    \widetilde R^\top & {\bf 0}
\end{bmatrix}P\right)$ is a NMF of the same rank.

Consider the singular value decomposition
\begin{equation} \label{svd} N \ge 0, \quad N = \sum_{i=1}^r  \bu_i \sigma_i \bv_i^\top.
\end{equation}
Then, it follows from the Perron-Frobenius theory that the dominant singular vectors $\bu_1$ and $\bv_1$ can be scaled to be nonnegative. Observe that, when $\sigma_1 = \sigma_2$, $\bu_1$ and $\bv_1$ are no longer uniquely defined (up to a sign); nonetheless, \cite[Theorem 2.1.1]{BermanP} guarantees that there exists a nonnegative choice of vectors in the dominant singular spaces. To simplify the notation, we define 
$ \hat\bu_i:=\bu_i\sigma_i^\frac12$ and $\hat\bv_i:=\bv_i\sigma_i^\frac12$ and rewrite this decomposition as follows
\begin{equation} \label{svd2}  N = \sum_{i=1}^r  \hat\bu_i \hat\bv_i^\top.
\end{equation} 
We point out that $\hat\bu_1$ and $\hat\bv_1$ are still nonnegative and nonzero vectors. In addition, Theorem \ref{thm:directsum} below discusses when $\hat \bu_1,\hat \bv_1$ are in fact both positive.

\begin{theorem}\label{thm:directsum}
    Let $0 \leq N \in \R^{m \times n}$, without any zero rows or columns, have a scaled singular value decomposition as in \eqref{svd2}, choosing without loss of generality $\hat \bu_1,\hat \bv_1 \geq 0$. Then, precisely one of the following is true:
    \begin{enumerate}
        \item $\hat \bu_1$ and $\hat \bv_1$ can be chosen to be both positive;
        \item $\hat \bu_1$ and $\hat \bv_1$ both necessarily have some zero entries.
    \end{enumerate}
    In addition, if $\hat \bu_1$ and $\hat \bv_1$ cannot be chosen positive, then $N$ is permutation equivalent to the direct sum of (at least) two non-empty nonnegative matrices.
\end{theorem}
\begin{proof}
Suppose first that there we can choose $\hat \bv_1 > 0$. Then, $\sigma_1 \hat \bu_1 = N \hat \bv_1 > 0$ (since $N \geq 0$ cannot have a zero row), implying that $\hat \bu_1 > 0$. An analogous argument, via $\sigma_1 \hat \bv_1 = N^\top \hat \bu_1$, shows the reverse implication. Hence, $\hat \bv_1$ can be chosen positive if and only if $\hat \bu_1$ can be chosen positive.

    Suppose now that $\hat \bv_1$ and $\hat \bu_1$ cannot be chosen positive, i.e., they necessarily both have some zero entries. By the Perron-Frobenius Theorem, the nonpositivity of $\hat \bv_1$ implies \cite[Theorem 2.1.3]{BermanP} that $N^\top N$ is reducible, and hence there exists a permutation matrix $P_c$ and non-empty nonnegative (and positive semidefinite) matrices $\widehat S,\widetilde S$ such that $P_c^\top N^\top N P_c = \widehat S \oplus \widetilde S$. Thus, the subspace generated by the leftmost columns of $NP_c$ is orthogonal to the subspace generated by the rightmost columns of $NP_c$. Taking into account that all these columns are nonnegative, we conclude that there must exist a second permutation matrix $P_r$ and non-empty nonnegative matrices $\widehat N,\widetilde N$  such that $P_rNP_c=\widehat N \oplus \widetilde N$ (and $\widehat S = \widehat N^\top \widehat N$, $\widetilde S = \widetilde N^\top \widetilde N$).
    \end{proof}

    \begin{remark}
        The last part of the statement of Theorem \ref{thm:directsum} gives a necessary, but not sufficient, condition for the impossibility of choosing positive dominant singular vectors. In other words, it is possible that $N=\widehat N \oplus \widetilde N$ but nonetheless $\hat \bu_1$ and $\hat \bv_1$ can be chosen positive. This can happen only if $\sigma_1=\|\widehat N\|_2=\|\widetilde N\|_2=\sigma_2$. For example, take $\widehat N=\widetilde N=1$ and $$\widehat U = \widehat V = \begin{bmatrix}
            \cos(\alpha) & \sin(\alpha)\\
            \sin(\alpha) & -\cos(\alpha)
        \end{bmatrix}, \qquad \alpha \in \left]0,\frac{\pi}{2}\right[.$$
        Then, clearly $0 \leq I_2 = \widehat N \oplus \widetilde N = \widehat U \widehat V^\top$ yet $\hat \bu_1 >0$ and $\hat \bv_1 > 0$.
    \end{remark}

    By repeated applications of Theorem \ref{thm:directsum}, we henceforth assume with no loss of generality that $\hat \bu_1,\hat \bv_1$ are both positive. Indeed, if $N=P_r^\top (\widehat N \oplus \widetilde N)P_c^\top$ for some permutation matrices $P_r,P_c$ and $\widehat N,\widetilde N \geq 0$, then we can compute NMFs $\widehat N=\widehat L \widehat R^\top$ and $\widetilde N=\widetilde L \widetilde R^\top$, and we observe that \[ N=\underbrace{\left(P_r^\top(\widehat L \oplus \widetilde L)\right)}_{:=L}\underbrace{\left((\widehat R^\top \oplus \widetilde R^\top)P_c^\top\right)}_{:=R^\top} \]
    is a NMF. Note that, for positive dominant singular vectors, the orthogonality of the pairs $(\bu_1,\bu_2)$ and $(\bv_1,\bv_2)$ guarantees that $\bu_2$ and $\bv_2$ 
must have both positive and negative components. In formulae,
\begin{equation}\label{eq:bothpandm}
    \exists \ i,j,k,\ell \ s.t. \ u_{i,2} > 0, \ u_{j,2} < 0, \ v_{k,2} > 0, \ v_{\ell,2} < 0. 
\end{equation} 
Clearly, if $r \geq 2$, then \eqref{eq:bothpandm} also applies to the scaled singular vectors $ \hat \bu_2, \hat \bv_2$.

\subsection{Case of nonnegative rank $2$}

We first assume that the rank of $N$ is $2$, i.e., \eqref{svd} holds with $r=2$ and $\sigma_1\geq \sigma_2 > 0$. To stress this rank property, below we denote such a rank-$2$ nonnegative matrix as $N_2$, while $M_2$ stands for a generic (possibly not nonnegative) matrix of rank $2$. Lemma \ref{lem:conditions} below studies a rank-$2$ matrix $M_2$ whose dominant singular vectors are both positive, and provides necessary and sufficient conditions for $M_2$ to be nonnegative.

\begin{lemma} \label{lem:conditions}
Consider the rank $2$ matrix $M_2$, without any zero row or column,
\begin{equation} \label{eq:svd2} M_2:= \widehat U\widehat V^\top := \left[\begin{array}{cc}\hat \bu_1 & \hat \bu_2\end{array}\right]\left[\begin{array}{cc}\hat \bv_1^\top \\ \hat \bv_2^\top\end{array}\right] 
\end{equation}
where $\hat \bu_1 > 0$, $\hat \bv_1 > 0$, and the pairs of vectors $(\hat \bu_1,\hat \bu_2)$ and $(\hat \bv_1,\hat \bv_2)$ are orthogonal.
Denote the $i$-th entry of $\hat \bu_k$ (resp. $\hat \bv_k$) by $\hat u_{i,k}$ (resp. $\hat v_{i,k}$) for $k=1,2$ and define
\begin{equation*} \minuh:=\min_i\frac{\hat u_{i,2}}{\hat u_{i,1}}, \;\;
\maxuh:=\max_i\frac{\hat u_{i,2}}{\hat u_{i,1}}, \; \;
\minvh:=\min_j\frac{\hat v_{j,2}}{\hat v_{j,1}},  \;\;  
\maxvh:=\max_j\frac{\hat v_{j,2}}{\hat v_{j,1}}.  
\end{equation*}
Then, $M_2$ is nonnegative if and only if
\begin{equation} \label{maxmin}  
\maxuh \minvh \ge -1, \quad \minuh\maxvh \ge -1. 
\end{equation}
\end{lemma}
\begin{proof}

The matrix $M_2$ is nonnegative if and only if
$$ \hat u_{i,1} \hat v_{j,1}+ \hat u_{i,2} \hat v_{j,2}\ge 0, \quad \forall i,j.
$$ 
Since the vectors $\hat \bu_1$ and $\hat \bv_1$ are nonnegative, this is equivalent to 
\begin{equation}\label{eq:claim}
    \left(\frac{\hat u_{i,2}}{\hat u_{i,1}}\right) \cdot \left(\frac{\hat v_{j,2}}{\hat v_{j,1}}\right)\ge -1, \quad \forall i,j.
\end{equation}

Each inequality in \eqref{eq:claim} is automatically satisfied either if the factors $\frac{\hat u_{i,2}}{\hat u_{i,1}}$ and $\frac{\hat v_{j,2}}{\hat v_{j,1}}$ have the same sign, or if one of them is zero. When they have opposite signs all inequalities will be satisfied if and only if the two extremal conditions \eqref{maxmin} are satisfied.
\end{proof}

The following general Theorem \ref{thm:vectors} is related to Lemma \ref{lem:conditions}.

\begin{theorem}\label{thm:vectors}
    Let ${\bf a},{\bf b} \in \R^n$ ($n \geq 2$) be orthonormal vectors, such that ${\bf a}$ does not have any zero entry. For $i=1,\dots,n$, let $q_i=b_i/a_i$. Then,
    \[ \max_i q_i (- \min_i q_i) \geq 1.  \]
\end{theorem}
\begin{proof}
Define $x:=-\min_i q_i$, $y:=\max_i q_i$, and note $x,y>0$. Denoting Dirac's measure centered on $r$ by $\delta_r$, observe that $\mu_Q:=\sum_{i=1}^n a_i^2 \delta_{q_i}$ is a probability measure (because $\sum_{i=1}^n a_i^2=1$) with support $\subseteq [-x,y]$. Furthermore, the associated random variable $Q$ has mean $\sum_{i=1}^n a_i^2 q_i=0$ and variance $\sum_{i=1}^n a_i^2 q_i^2 = 1$.

Consider now the set $\mathcal{S}$ of probability measures with support $\subseteq [-x,y]$ and zero mean, so that $\mu_Q \in \mathcal{S}$, and let $X$ be the random variable whose distribution $\mu_X$ has maximal variance within $\mathcal{S}$. Since $\mathcal{S}$ is a convex set and the second moment is a linear function of the measure, $\mu_X$ is an extreme point of $\mathcal{S}$ \cite{Rockafellar}, i.e., a measure that cannot be written as a convex combination of other two measures in $\mathcal{S}$. These are precisely the $2$-atomic distributions \cite{Winkler}. It is then easy to see that 
\[  \mu_X = \frac{y}{x+y} \delta_{-x} + \frac{x}{x+y} \delta_y.\] 
Therefore,
\[  1=\mathbb{E}[Q^2] \leq \mathbb{E}[X^2] = \frac{x^2y+xy^2}{x+y} = xy. \]
\end{proof}

By Lemma \ref{lem:conditions}, a necessary condition for nonnegativity of $M_2$ is $$\maxuh \maxvh \minuh \minvh \leq 1.$$ By scaling the vector ${\bf q}$ in the statement of Theorem \ref{thm:vectors}, we obtain instead a lower bound for this product:
\begin{equation}\label{eq:lowerbound}
    \min\left\{ (-\minuh \maxuh),(-\minvh\maxvh)   \right\}     \geq \frac{\sigma_2}{\sigma_1}  \Rightarrow \maxuh \maxvh \minuh \minvh \geq \left(\frac{\sigma_2}{\sigma_1}\right)^2.
\end{equation}   
As a consequence, when $\sigma_1=\sigma_2$ and $M_2$ is nonnegative, \eqref{maxmin} must in fact hold with equality.

\subsection{A geometric reformulation}
We now reinterpret Lemma \ref{lem:conditions} geometrically by using angular parameters. Let us define 
 $$
 \psi_i:=\arctan\left(\frac{\hat u_{i,2}}{\hat u_{i,1}}\right), \qquad 
 -\pi/2 < \psi_i < \pi/2 
 $$
 and
 $$
 \phi_j:=\arctan\left(\frac{\hat v_{j,2}}{\hat v_{j,1}}\right), \qquad 
 -\pi/2 < \phi_j < \pi/2.
 $$
The condition \eqref{maxmin} then expresses that 
\begin{equation*}    |\psi_i-\phi_j|\le \pi/2 \quad \forall i,j \quad \Leftrightarrow \quad \hat u_{i,1}\hat v_{j,1} +\hat u_{i,2}\hat v_{j,2} \ge 0 \quad \forall i,j.
\end{equation*} 

When $\hat u_{i,2}$ and $\hat v_{j,2}$ have the same sign or one of them is zero, the condition is automatically satisfied since the vectors $[\cos(\psi_i),\sin(\psi_i)]^\top$ and $[\cos(\phi_j),\sin(\phi_j)]^\top$ are in the same closed quadrant. When $\hat u_{i,2}$ and $\hat v_{j,2}$ have different signs, 
$[\cos(\psi_i),\sin(\psi_i)]^\top$ and $[\cos(\phi_j),\sin(\phi_j)]^\top$
are in adjacent quadrants. The matrix $M_2$ is then nonnegative if and only if 
$ |\psi_i-\phi_j|\le \pi/2$ for all $i,j$. 
Moreover, if we define the positive scaling factors
$d_i^{(u)}:=\sqrt{\hat u_{i,1}^2+\hat u_{i,2}^2}$ and
$d_j^{(v)}:=\sqrt{\hat v_{j,1}^2+\hat v_{j,2}^2}$, and the 
diagonal matrices 
$$ D_u:=\mathrm{diag}\left(d_1^{(u)},\ldots,d_m^{(u)}\right), \qquad 
 D_v:=\mathrm{diag}\left(d_1^{(v)},\ldots,d_n^{(v)}\right),
$$
then we can write $M_2$ as follows
\begin{equation} \label{eq:DcosD}
 M_2:= D_u \left[\begin{array}{ccc}
 \cos(\psi_1-\phi_1) & \ldots & \cos(\psi_1-\phi_n)\\
 \vdots & \ddots & \vdots \\ \cos(\psi_m-\phi_1) & \ldots & \cos(\psi_m-\phi_n)  \end{array} \right]D_v.
\end{equation}
Observe that, once more, the fact that the scaling factors are positive follows from the assumption that $M_2$ does not have any zero row or zero column.

\begin{remark}  \label{rem:flip}
Note that in the factorization \eqref{eq:svd2} we can change the sign of the vectors $\hat\bu_2$ and $\hat\bv_2$ without affecting the product $M_2$. This would flip the signs of all the angles, $\phi_i \leftarrow - \phi_i, \psi_j \leftarrow -\psi_j$, without affecting the constraints in Lemma \ref{lem:conditions}.
\end{remark}

Furthermore, the angles $\psi_i$ and $\phi_j$ in \eqref{eq:DcosD} can be ordered non-decreasingly by introducing permutation matrices $P_\psi M_2 P_\phi$ on rows and columns of $M_2$, that obviously do not affect nonnegativity of the factorization. The angles will be ordered non-decreasingly if and only if the ratios $\hat u_{i,2}/\hat u_{i,1}$ and $\hat v_{i,2}/\hat v_{i,1}$ are ordered non-decreasingly, because
the function $\arctan$ is strictly monotone increasing on $\R$.
After this permutation, the matrix $P_\psi M_2 P_\phi$ will be nonnegative if and only if the extreme elements
$\cos(\psi_1-\phi_n)$ and $\cos(\psi_m-\phi_1)$ are nonnegative.
Moreover, if $P_\psi M_2 P_\phi$ has negative elements, then they must be grouped around the elements $(1,n)$ and $(m,1)$ and form a staircase-type pattern around these corners. 
This observation is also independent of the possible flipping of the signs of the angles as described in Remark \ref{rem:flip}.
 
Without loss of generality, we could therefore assume that $P_\psi$ and $P_\phi$ are identity matrices, i.e., that the matrix $M_2$ has reordered rows and columns.
The sign pattern of $M_2$ in \eqref{eq:DcosD} is then that of the middle factor since $D_u \succ 0, D_v \succ 0$ and it  can then be partitioned into four blocks 
$$   M_2 = 
\left[ \begin{array}{cc} M_{11} & M_{12} \\
M_{21} & M_{22} \end{array}\right]
$$
by grouping the elements based on which quadrants $\psi_i$ and $\phi_j$ belong: the diagonal blocks (which are possibly non-square) are positive and correspond to the 
angle pairs $(\psi_i,\phi_j)$ where both angles lay in $]-\pi/2,0[$ (for $M_{11}$) or both in $[0,\pi/2[$ (for
$M_{22}$). The blocks $M_{21}$ and $M_{12}$ contain the elements where the angles lay in different quadrants. The sign pattern of 
$M_2$ is then described in the following proposition.

\begin{prop}\label{prop:negcorners}
    Let $M_2$ be a rank-2 matrix, without any zero row or column and with positive dominant singular vectors, with a factorization of form \eqref{eq:DcosD} where $$\pi/2 > \psi_1 \geq \dots \geq \psi_m > -\pi/2 \quad and \quad \pi/2 > \phi_1 \geq \dots \geq \phi_n > -\pi/2.$$ 
    Then the negative entries of $M_2$ are grouped in the upper-right and lower-left corners of the matrix, and satisfy a staircase-like pattern, that is
    \[ \begin{cases}
        M_{12}(i,j) < 0 \implies M_{12}(k,\ell) < 0,~ \forall  k \leq i, \ell \geq j \\
        M_{21}(i,j) < 0 \implies M_{21}(k,\ell) < 0,~ \forall  k \geq i, \ell \leq j. 
    \end{cases}  \]   
\end{prop}
\begin{proof}
    It is enough to inspect the nonnegativity of the middle factor. An element of it is negative when
    \[ \cos(\psi_i - \phi_j) < 0 \Leftrightarrow |\psi_i - \phi_j| > \pi/2. \]
    This can only happen when the angles $\psi_i$ and $\phi_j$ lay in different quadrants, and so we are in one of the submatrices $M_{12}$ or $M_{21}$. Suppose we are in $M_{12}$, and so $\psi_i > \phi_j + \pi/2$. Since the angles are ordered, we have
    \[ \psi_k \geq \psi_i > \phi_j + \pi/2 \geq \phi_\ell + \pi/2, ~\forall k \geq i, \ell \leq j. \]
    A similar argument can be made when we are in $M_{21}$.
\end{proof}

\begin{prop}
    Let $M_2$ be as in Proposition \ref{prop:negcorners}. If $M_2$ is nonnegative, then it is of the form
    \[ M_2 = \begin{bmatrix}
        M_{11} & M_{12} & 0 \\
        M_{21} & M_{22} & M_{23} \\
        0 & M_{32} & M_{33}
    \end{bmatrix} \]
    where the nonzero blocks $M_{ij}$ are positive.
\end{prop}
\begin{proof}
    As the entries are ordered as in Proposition \ref{prop:negcorners}, the smallest entries are packed to the top-right and bottom-left corners. The zero entries must form a block of zeros, because the matrix is of rank $2$. Indeed, suppose that in the bottom-left corner the zero entries do not form a block. Then there would be a $3\times 3$ submatrix of the form (the lines split elements of the bottom-left block of $M_2$ from the rest)
    \[\begin{array}{cc|c}
        \star & \star & \star\\
        \hline
        0&\star&\star\\\
        0&0&\star
    \end{array},\] 
where the symbol $\star$ denotes a positive real number whose precise value is irrelevant. But this implies $\rank ~M_2 \geq 3$, which is a contradiction. An analogous argument holds for the top-right corner. Using also that $M_2$ cannot have zero rows and columns, the statement follows.
\end{proof}

\subsection{A parametrization of rank 2 NMFs}

When $M_2 \geq 0$, consider now a nonnegative rank-two factorization $M_2=L R^\top$ with $L,R \geq 0$, whose existence is well estabilished \cite{Gillis}. The factors $L$ and $R$ are then related to the factors $\widehat U$ and $\widehat V$ of Lemma \ref{lem:conditions} by an invertible transformation $X \in \R^{2 \times 2}$~:
$$   LR^\top=\left(\left[\begin{array}{cc}\hat \bu_1 & \hat \bu_2\end{array}\right]X\right)\left(X^{-1}\left[\begin{array}{cc}\hat \bv_1^\top \\ \hat \bv_2^\top\end{array}\right]\right), \quad \det X \neq 0.
$$
Note further that, by \eqref{eq:bothpandm} (recall that we assume $\hat \bu_1 > 0, \hat \bv_1 >0$), another necessary condition for the nonnegativity of $L$ and $R$ is that $X$ and $X^{-1}$ cannot contain zero entries. Theorem \ref{th:main} below shows that $X$ can always be parameterized using four positive elements in a convenient factorization $X=TDP$ where $D$ is diagonal with positive diagonal elements and $P$ is a permutation matrix. 
\begin{theorem} \label{th:main}
Let $N_2\ge 0$ have rank $2$, not have any zero row or column, and be given by its scaled singular value decomposition \eqref{eq:svd2}. Assume moreover that $\hat \bu_1, \hat \bv_1$ are both positive. Then the factorization
\begin{equation} \label{eq:LRUV} L R^\top=\left(\left[\begin{array}{cc}\hat \bu_1 & \hat \bu_2\end{array}\right]TDP\right)\left((TDP)^{-1}\left[\begin{array}{cc}\hat \bv_1^\top \\ \hat \bv_2^\top\end{array}\right]\right)
\end{equation}
with $D=\begin{bmatrix}
    d_1 & 0\\
    0 & d_2
\end{bmatrix}, \;d_1,d_2>0$, $T=\left[\begin{array}{cc} 1 & t_2 \\ t_1 & -1 \end{array}\right], \; t_1,t_2>0$, and $P$ a $2\times 2$ permutation matrix, describes the whole family of nonnegative rank-two factorizations of $N_2$, provided
$$   \maxvh \le t_1  \le  -1/\minuh, $$
$$    \maxuh \le t_2  \le  \: -1/\minvh. $$
\end{theorem}
\begin{proof}
Since $X$ cannot have any zero entry, we can always factorize it as $X=TDP$ where $P$ is a $2\times 2$ permutation, and
$$   T=\left[\begin{array}{cc} e_1 & t_2 \\ t_1 & e_2 \end{array}\right], \quad  D=\left[\begin{array}{cc} d_1 & 0 \\ 0 & d_2  \end{array}\right], \quad d_1,d_2>0, \quad |e_1|=|e_2|=1.$$
If this is the case then
$$
T^{-1}=\left[\begin{array}{cc} -e_2 & t_2 \\ t_1 & -e_1 \end{array}\right] \cdot \frac{1}{(t_1t_2-e_1e_2)}.
$$
In addition, we can assume without loss of generality that $t_1t_2-e_1e_2$ (which is minus the determinant of $T$)
is positive, for if $\det T >0$ then we can simply change the choice of the permutation $P$ in the product $TDP$.

Let us now impose that the products 
$$
\left[\begin{array}{cc}\hat \bu_1 & \hat \bu_2\end{array}\right]T \quad \mathrm{and}\quad 
\left[\begin{array}{cc}\hat \bv_1 & \hat \bv_2\end{array}\right]T^{-\top}
$$
are nonnegative. Here, we can dismiss the scaling factor $D$ and permutation matrix $P$, because they do not affect the nonnegativity of the factors. The conditions for $L\ge 0$ and $R\ge 0$ then become
\begin{equation} \label{eq:L}  L_{i,1}= \hat u_{i,1}e_1+\hat u_{i,2} t_1 \ge 0,  \quad   L_{i,2}=\hat u_{i,1}t_2 +\hat u_{i,2}e_2  \ge 0,   \quad \forall i \in [m],
\end{equation}
\begin{equation} \label{eq:R}   R_{j,1}= -\hat v_{j,1}e_2+\hat v_{j,2} t_2 \ge 0, \quad   R_{j,2}= \hat v_{j,1}t_1 -\hat v_{j2}e_1  \ge 0,    \quad \forall j \in [n].
\end{equation}
Combining \eqref{eq:L}, \eqref{eq:R} and \eqref{eq:bothpandm}, it follows that $e_1 >0$, $t_1 >0$, $-e_2>0$, $t_2 >0$. This implies $e_1=1$, $e_2=-1$ and
\begin{eqnarray} 0<  \maxvh \le  & t_1 &  \le -1/\minuh, \label{eq:t1} \\
  0< \maxuh \le & t_2 & \le  -1/\minvh.  \label{eq:t2}
\end{eqnarray}
It follows from Lemma \ref{lem:conditions} that these conditions are compatible with the nonnegativity of $N_2$.
\end{proof}

\begin{corollary}  \label{cor:unique}
The nonnegative factorization $N_2=L R^\top=(\widehat UT)(T^{-1}\widehat V^\top)$ given in \eqref{eq:LRUV} is unique (up to a diagonal scaling $D$ and a permutation $P$) if and only if 
\begin{equation} \label{diag}
 \maxuh \minvh = -1, \quad \maxvh\minvh=-1.
\end{equation} 
\end{corollary}
\begin{proof}
Equations \eqref{diag} follow from Theorem \ref{th:main} and Lemma \ref{lem:conditions}, by focusing on the case where the allowed intervals for $t_1$ and $t_2$ reduce to single points.
\end{proof}
Corollary \ref{cor:unique} corresponds to \cite[Theorem 4.2]{Gillis}, since the
condition is equivalent to the existence of a diagonal submatrix in the rows and columns of $N_2$ corresponding to the two equalities in \eqref{diag}. We note that our proof is elementary and relies on studying a number of equalities/inequalities; the proof in \cite{Gillis} invokes instead more advanced geometric tools such as convex hulls. In addition, we note that by \eqref{maxmin} and \eqref{eq:lowerbound} the uniqueness conditions are always satisfied when $\sigma_1= \sigma_2$.

Lemma \ref{lem:alpha} revisits Theorem \ref{th:main} by combining it with the decomposition 
\eqref{eq:DcosD}, and giving conditions for when the corresponding rank-two matrix is nonnegative.

\begin{lemma} \label{lem:alpha}
Consider the rank-two matrix $M_2:=\left[\begin{array}{cc}
\hat \bu_1 & \hat \bu_2 \end{array}\right] \left[\begin{array}{cc}
\hat \bv_1 & \hat \bv_2 \end{array}\right]^\top$, without any zero row or column and 
where the vectors $\hat \bu_1$ and  $\hat \bv_1$ are positive, and the vectors $\hat \bu_2$ and  $\hat \bv_2$ both have elements of different signs.
Define the vectors 
$$ {\bf {\bf \Psi}}:=\left[\begin{array}{c}
\psi_1 \\ \vdots \\ \psi_m \end{array}\right] , \quad {\bf \Phi}:=\left[\begin{array}{c}
\phi_1 \\ \vdots \\ \phi_n \end{array}\right], \quad \1_k :=
\left[\begin{array}{c} 1 \\ \vdots \\ 1 \end{array}\right] \in \R^k
$$
where $$\tan(\psi_i)=\frac{\hat u_{i,2}}{\hat u_{i,1}}, \quad 
\tan(\phi_j)=\frac{\hat v_{j,2}}{\hat v_{j,1}}, \quad -\pi/2< \psi_i,\phi_j < \pi/2, $$
and define
$D_u:=\mathrm{diag}\left(d_1^{(u)},\ldots,d_m^{(u)}\right), \; 
 D_v:=\mathrm{diag}\left(d_1^{(v)},\ldots,d_n^{(v)}\right),$
where
$$ d_i^{(u)}:=\sqrt{\hat u_{i,1}^2+\hat u_{i,2}^2}, \quad
d_j^{(v)}:=\sqrt{\hat v_{j,1}^2+\hat v_{j,2}^2}.$$ 
Then $M_2$ can be written as
\begin{equation} \label{eq:DPhiPsi}
M_2 =  D_u [\cos({\bf \Psi}),\sin({\bf \Psi})] [\cos({\bf \Phi}),\sin({\bf \Phi})]^\top D_v  =  D_u \cos([{\bf \Psi}\1_n^\top - \1_m {\bf \Phi}^\top]) D_v.
\end{equation}
Moreover, the matrix 
$M_2$ is nonnegative if and only if there exist angles $\theta_1$ and $\theta_2$ in $]0,\pi/2[$ such that 
\begin{equation}
    \label{eq:theta12}
(\theta_1-\pi/2)\1_m \le {\bf \Psi} \le \theta_2 \1_m , \quad 
(\theta_2-\pi/2)\1_n \le {\bf \Phi} \le \theta_1 \1_n. 
\end{equation} 
If this is the case, then a nonnegative factorization of the middle factor
$\cos([{\bf \Psi}\1_n^\top - \1_m {\bf \Phi}^\top])$ is given by
\begin{equation} \label{eq:factors}
 \frac{1}{\cos(\alpha_2-\alpha_1)}.\left[\cos(\alpha_1\1-{\bf \Psi}), \sin(\alpha_2\1-{\bf \Psi})\right]
 \left[\cos(\alpha_2 \1 -{\bf \Phi}),\sin(\alpha_1\1-{\bf \Phi})\right]^\top,
\end{equation}
for any parameters $\alpha_1$ and $\alpha_2$ in the intervals
\begin{equation}
\label{eq:int}
\max_j \phi_{j} \le \alpha_1 \le \min_i \psi_{i}+\pi/2, \quad
 \max_i \psi_i \le \alpha_2 \le \min_j \phi_{j}+\pi/2,
\end{equation}
\end{lemma}
\begin{proof}
The conditions \eqref{eq:t1} and \eqref{eq:t2} translate into
the inequalities \eqref{eq:int},
which are all quantities in the open interval $]0,\pi/2[$. 
We can then define the angles $\theta_1$, $\theta_2$ in
$]0,\pi/2[$ as follows
$$   \theta_1:= \max_j \phi_j \le \min_i \psi_i+\pi/2, \quad 
\mathrm{and} \quad \theta_2:=\max_i \psi_i \le \min_j \phi_j+ \pi/2,
$$
to satisfy \eqref{eq:theta12}. The updating transformation $(TD)(TD)^{-1}$ of Theorem \ref{th:main}
is replaced here by an equivalent factorization of $I_2$
$$  I_2=\frac{1}{\cos(\alpha_2-\alpha_1)} \left[\begin{array}{cc} \cos(\alpha_1) & \sin(\alpha_2) \\ \sin(\alpha_1) & -\cos(\alpha_2) \end{array}\right]
 \left[\begin{array}{cc} \cos(\alpha_2) & \sin(\alpha_2) \\ \sin(\alpha_1) & -\cos(\alpha_1) \end{array}\right],
$$
which yields the desired rank-two factorization. Moreover, the factors are nonnegative since $\alpha_2-\psi_i$ and $\alpha_1-\phi_j$ are in $]0,\pi/2[$,
and $\alpha_2-\alpha_1$ lies in $]-\pi/2,\pi/2[$.
\end{proof}

We note in passing that in Lemma \ref{lem:alpha} we could have chosen other values for $\alpha_1$ and $\alpha_2$,
for instance the midpoints of the intervals \eqref{eq:int}~:
$$
    \widehat\alpha_1:= \frac12(\min_i \psi_{i} + \max_j \phi_{j}+\pi/2),\quad \widehat\alpha_2:=
     \frac12(\min_j \phi_j + \max_i \psi_{i}+\pi/2).
$$
Such a choice guarantees that the new angles are not on the boundary of the intervals \eqref{eq:int} when they have a non-vanishing interior.

\begin{remark} \label{rem:balanced}
The factorization $M_2=L R^\top$ in \eqref{eq:DPhiPsi} implies that the Frobenius norms of the left and right factors $L$ and $R$ are equal to $\|D_u\|_F$ and $\|D_v\|_F$, respectively. Moreover, the initial factorization that is based on the singular value decomposition of $N_2$ implies that
$$ \|D_u\|_F = \|D_v\|_F = \sqrt{\sigma_1+\sigma_2},$$
which means that the representation is in a sense ``balanced".
We will see that during the updating of these factorizations, these norms can both only decrease and that they remain approximately balanced. 
\end{remark}

\subsection{Special rank-$2$ nonnegative factorizations}

The special types of rank-$2$ factorizations described in this subsection have a direct application in scaled stochastic block modeling with two roles,
 where the expected value of the weighted adjacency matrix has the form \cite{Barbarino,Marchand}
\[   \mathbb{E}[N] = \left[\begin{array}{cc} {\bf d}_{\ell 1} & {\bf 0}_{n_1} \\ {\bf 0}_{n_2} & {\bf d}_{\ell 2}\end{array}\right] M \left[\begin{array}{cc} {\bf d}_{r1}^\top & {\bf 0}^\top_{n_2} \\ {\bf 0}^\top_{n_1} & {\bf d}_{r2}^\top \end{array}\right] ,\]
where ${\bf d}_{\ell i} \geq 0 \in \mathbb{R}^{n_i}$, ${\bf d}_{r i} \geq 0 \in \mathbb{R}^{n_i}$, $M \geq 0 \in 
\mathbb{R}^{2\times 2}. $

Such rank-$2$ factorizations of a nonnegative matrix can also be seen as a ``three way factorization" $N_2 := L M R^\top$ using a nonnegative middle matrix as well as nonnegative left and right factors. This possibility is discussed, for example, in \cite{DacheVG24,Gillis}. In such a case, one may characterize further properties that are desired in the three factors (besides nonnegativity), since there are additional degrees of freedom with respect to a traditional NMF. In the case of scaled stochastic block modeling \cite{Barbarino,Marchand}, one is looking for a factorization
$$ N_2= L M R^\top , \quad  \mathrm{diag}(L^\top L)=\mathrm{diag}(R^\top R)= I_2 $$
and where $\|L^\top L-I_2\|_2$ and $\|R^\top R-I_2\|_2$ are minimized. Informally, this means that the columns of $L$ and $R$ are ``as close as possible" to being orthonormal. We show in this section that such a factorization always exists and is easy to compute. As another example, if one is given a nonnegative matrix $0 \leq N \in \R^{n \times n}$ that is symmetric and positive semidefinite, one often wishes to approximate it by another rank-$2$ matrix with the same properties. Semidefinite matrices that admit a 
nonnegative factorization $N=LL^\top(=L I L^\top)$, where $L\ge 0$, are said to be {\em completely positive}. They are discussed, e.g., in \cite{BermanP,Diananda,GW}, and studied in \cite{KG12} using ideas similar to the ones proposed here. When $n \leq 4$, a matrix is completely positive if and only if it is a nonnegative positive semidefinite matrix; but if $n \geq 5$, counterexamples exist \cite{Diananda,GW}.  

Theorem \ref{th:threeway} shows how, for the rank-$2$ case, our analysis extends to general three way factorizations.

\begin{theorem} \label{th:threeway}
Let $N_2\ge 0$ have rank $2$, not have any zero row or column, and be given by its scaled singular value decomposition \eqref{eq:svd2}, where the vectors $\hat\bu_1$
and $\hat\bv_1$ are positive. Then the factorization
\begin{equation} \label{eq:LR} N_2= L M R^\top=\left(\left[\begin{array}{cc}\hat \bu_1 & \hat \bu_2\end{array}\right]T_LD_L\right)\left(D_RT_R ^\top T_LD_L \right)^{-1} \left(D_RT_R^\top\left[\begin{array}{cc}\hat \bv_1^\top \\ \hat \bv_2^\top\end{array}\right]\right)
\end{equation}
with $T_L=\left[\begin{array}{cc} 1 & \underline{t}_2 \\ \overline{t}_1 & -1 \end{array}\right]/\sqrt{\overline{t}_1\underline{t}_2+1},$ and $T_R^\top=\left[\begin{array}{cc} 1 & \overline{t}_2 \\ \underline{t}_1 & -1  \end{array}\right]/\sqrt{\underline{t}_1\overline{t}_2+1},$ is a three way rank two and nonnegative factorization of $N_2$, for all parameters $\underline{t}_1, \overline{t}_1, \underline{t}_2, \overline{t}_2$ satisfying
\begin{eqnarray}
    \label{eq:tangents}
0 < \maxuh \le \underline{t}_1 \le \overline{t}_1\le -1/\minvh, \\ \nonumber
0 <  \maxvh \le \underline{t}_2 \le \overline{t}_2 \le -1/\minuh . 
\end{eqnarray}
and where the matrices $D_L$ and $D_R$ are positive diagonal scaling matrices.
\end{theorem}
\begin{proof}
The transformations $T_L$ and $T_R$ follow the rules of Theorem \ref{th:main}, but have been scaled so that $T_L^{-1}=-T_L$
and  $T_R^{-1}=-T_R$.
It remains to prove that the middle matrix $M:=(D_RT_R^\top T_LD_L)^{-1}$ is nonnegative. That follows from the equalities  $$ D_L M D_R = T_L^{-1} T_R^{-\top} = T_L T_R^\top  $$
and 
$$ \left[\begin{array}{cc} 1 & \underline{t}_2 \\ \overline{t}_1 & -1 \end{array}\right] \left[\begin{array}{cc} 1 & \overline{t}_2 \\ \underline{t}_1 & -1  \end{array}\right] =
\left[\begin{array}{cc} 1 +\underline{t}_2\underline{t}_1& \overline{t}_2 -\underline{t}_2 \\ \overline{t}_1-\underline{t}_1 & 1+\overline{t}_1\overline{t}_2  \end{array}\right] \ge 0.
$$
\end{proof}

\begin{remark}
It follows from the construction of the $LMR^\top$ factorization \eqref{eq:LR} given in Theorem \ref{th:threeway} that this describes the complete parameterization of such a factorization, except for additional permutations of the form $(LP_L)(P_L^\top M P_R)(P_R^\top R^\top).$
\end{remark}

Let us now look at the angles $\theta_L$ and $\theta_R$ between the two columns $\ell_i$ of $L$ and $r_i$ of $R$, defined via~:
$$ \cos^2\theta_L:= (\ell_1^\top\ell_2)^2/[(\ell_1^\top\ell_2)(\ell_2^\top\ell_2)], \quad \cos^2\theta_L:= (r_1^\top r_2)^2/[(r_1^\top r_2)(r_2^\top r_2)],
$$
If we introduce additional $2\times 2$ diagonal scalings $\widehat D_L\succ 0$ and $\widehat D_R\succ 0$ in the factorization $\widehat L \widehat M \widehat R^\top:=(L\widehat D_L)(\widehat D_L^{-1}M\widehat D_R^{-1})(\widehat D_RR^T)$ such that the scaled matrices 
$$ \widehat L^\top \widehat L :=\widehat D_L L^\top L \widehat D_L  \quad \mathrm{and} \quad  
\widehat R^\top \widehat R :=\widehat D_R R^\top R \widehat D_R$$
have unit diagonal elements\footnote{Note that the factor $\widehat D_L$ (resp. $\widehat D_R$) can be absorbed in $D_L$ (resp. $D_R$).} then the angles $\theta_L$ and $\theta_R$ can be linked to the so-called {\em orthogonality defects} of $\widehat L$ and $\widehat R$~: 
$$\mathrm{def}_L:=\|\widehat L^\top \widehat L-I_2\|_2=|\cos\theta_L|, \quad \mathrm{def}_R:=\|\widehat R^\top \widehat R-I_2\|_2=|\cos\theta_R|.$$
These can be expressed in function of the parameters 
$\underline{t}_1$, $\overline{t}_1$,
$\underline{t}_2$ and $\overline{t}_2$ as follows.
Let
$$ L^\top L=D_LT_L^\top\Sigma T_L D_L = \left[\begin{array}{cc}
\sigma_1 +\sigma_2 \overline{t}_1^2 & \sigma_1 \underline{t}_2 -\sigma_2 \overline{t}_1 \\ \sigma_1 \underline{t}_2 -\sigma_2 \overline{t}_1 & \sigma_2 +\sigma_1 \underline{t}_2^2
\end{array}\right],
$$
which implies that 
$$\cos^2\theta_L :=
(\sigma_1 \underline{t}_2 -\sigma_2 \overline{t}_1)^2/[(\sigma_1 +\sigma_2 \overline{t}_1^2)(\sigma_2 +\sigma_1 \underline{t}_2^2)]
$$
and $$\cos^2\theta_R :=
(\sigma_1 \underline{t}_1 -\sigma_2 \overline{t}_2)^2/[(\sigma_1 +\sigma_2 \overline{t}_2^2)(\sigma_2 +\sigma_1 \underline{t}_1^2)].
$$

In Proposition \ref{th:minimal}, we show how to minimize the orthogonality defects. To better interpret the statement of Proposition \ref{th:minimal}, recall that  by \eqref{eq:lowerbound} it always must hold $\sigma_1/\sigma_2 \geq -1/(\minuh \maxuh)$ and $\sigma_1/\sigma_2 \geq -1/(\minvh \maxvh)$.

\begin{prop}
 \label{th:minimal}
 In the notation and setting defined above, the functions $\cos^2\theta_L$
 and $\cos^2\theta_R$ reach their respective minima over the sets \eqref{eq:tangents} in the corner points, resp.,
$(\overline{t}_1=-1/\minvh$, $\underline{t}_2=\maxvh$), and
$(\overline{t}_2=-1/\minuh, \underline{t}_1=\maxuh)$.
Moreover, a zero defect of $\widehat L$ can only be obtained if $\sigma_1/\sigma_2=-1/(\minvh\maxvh)$ and a zero defect of $\widehat R$ can only be obtained if $\sigma_1/\sigma_2=-1/(\minuh\maxuh)$.
\end{prop}
\begin{proof} 
We show how to minimize $\mathrm{def}_L$; the argument for $\mathrm{def}_R$ is analogous. Since $L$ is a nonnegative matrix, or alternatively by combining \eqref{eq:lowerbound} and \eqref{eq:tangents}, we have $\sigma_1 \underline{t}_2 - \sigma_2 \overline{t}_1 \geq 0$. On the other hand, by a direct computation, we get that when $\cos^2\theta_L=0$ its partial derivatives are also zero, and otherwise

\[  \mathrm{sign} \left(\frac{\partial \cos^2 \theta_L}{\partial \underline{t}_2}\right) = \mathrm{sign} (\sigma_1 \underline{t}_2 - \sigma_2 \overline{t}_1) = - \mathrm{sign} \left(\frac{\partial \cos^2 \theta_L}{\partial \overline{t}_1}\right).         \]

Therefore, $\cos^2 \theta_L$ is minimized when $\overline{t_1}=-1/\minvh$ is at the right endpoint and $\underline{t}_2=\maxvh$ is at the left endpoint of their allowed intervals \eqref{eq:tangents}. Furthermore, the corresponding minimal defect is zero precisely when $\sigma_1 \underline{t}_2-\sigma_2 \overline{t_1}=\sigma_1 \maxvh + \sigma_2/\minvh=0$.\end{proof}

\begin{example} Consider the symmetric nonnegative and positive semidefinite matrix (given in factorized form)
$$  N_2 = LML^\top :=
 \left[\begin{array}{cr} \frac{4}{\sqrt{50}} & \frac12 \\ \frac{3}{\sqrt{50}} & \frac12 \\ \frac{3}{\sqrt{50}} & -\frac12 \\ \frac{4}{\sqrt{50}} & -\frac12 \end{array} \right]   \left[\begin{array}{cr} 6 & 2 \\ 2 & 3 \end{array}\right]
 \left[\begin{array}{ccrr} \frac{4}{\sqrt{50}} & \frac{3}{\sqrt{50}} & \frac{3}{\sqrt{50}} & \frac{4}{\sqrt{50}} \\
  \frac12 & \frac12 & -\frac12 & -\frac12 \end{array} \right],
$$
where $L^\top L=I_2$.
Its eigenvalues are $\lambda_1=7$ and $\lambda_2=2$ and its corresponding scaled eigenvectors and ratios, are
$$ \hat \bu_1 = 
L  \left[\begin{array}{cr} 2 \\ 1 \end{array}\right]\sqrt{7/5}, \quad \hat \bu_2 = L \left[\begin{array}{r} -1 \\ 2 \end{array}\right]\sqrt{2/5}, \quad \hat \bu_2./\hat \bu_1 \approx \left[\begin{array}{r} 0.1423 \\ 0.2282 \\ -2.1843 \\ -1.3255 \end{array}\right].$$
The above decomposition is not nonnegative, but any values of $\underline{t}$ and  $\overline{t}$ satisfying 
$$ 0.2282 \le \underline{t}\le \overline{t} \le 1/2.1843=0.4578
$$
will yield a three way nonnegative factorization. The one corresponding to the extremal values $\underline{t} = 0.2282$ and $\overline{t} = 0.4578$ yields the optimal nonnegative solution since $\sigma_1/\sigma_2=3.5 >  \overline{t}/\underline{t}$. The corresponding factors are 
$$ \widehat L \approx \left[\begin{array}{cr} 0.7548 & 0.1078 \\
 0.6470 & 0 \\  0 &  0.6470 \\ 0.1078 &  0.7548 \end{array}\right], \quad \widehat M \approx  \left[\begin{array}{cr} 6.3987 & 0.7883 \\
 0.7883 &2.3446 \end{array}\right]
$$
and the orthogonality defect is $\|\widehat L^\top \widehat L-I_2\|_2 \approx 0.1628.$ Note that each column of $\widehat L$ must contain at least one boundary element since otherwise one could 
rotate these columns further.
\end{example}

When we specialize Theorem \ref{th:threeway} to a symmetric nonnegative matrix $N_2$, we obtain the following theorem.
\begin{theorem} \label{th:threewsym}
Let the symmetric matrix $N_2\ge 0$ have rank $2$, not have any zero row or column, and be given by its eigenvalue decomposition 
$$ N_2= \left[\begin{array}{cc}\bu_1 & \bu_2\end{array}\right] \Lambda \left[\begin{array}{c} \bu_1^\top \\ \bu_2^\top \end{array}\right]= \left[\begin{array}{cc}\hat\bu_1 & \hat\bu_2\end{array}\right] S \left[\begin{array}{c} \hat\bu_1^\top \\ \hat\bu_2^\top \end{array}\right],
$$
where $\hat\bu_1=\bu_1\lambda_1^\frac12$, $\hat\bu_2=\bu_2|\lambda_2|^\frac12$, $\hat\bu_1$ is positive, and either $S=\begin{bmatrix}
    1&0\\
    0&1
\end{bmatrix}$ in the  semidefinite case or  $S=\begin{bmatrix}
        1&0\\
        0&-1
    \end{bmatrix}$ in the indefinite case. 
    Then the factorization
\begin{equation*} L M L^\top=\left(\left[\begin{array}{cc}\hat \bu_1 & \hat \bu_2\end{array}\right]TD\right)\left(DT ^\top S TD \right)^{-1} \left(DT^\top\left[\begin{array}{cc}\hat \bu_1^\top \\ \hat \bu_2^\top\end{array}\right]\right),
\end{equation*}
with $T=\left[\begin{array}{cc} 1 & \underline{t} \\ \overline{t} & -1 \end{array}\right]/\sqrt{\overline{t}\underline{t}+1}$ and $D \succ 0$ is a scaling diagonal matrix, 
 is a three way rank-$2$ nonnegative factorizations of $N_2$, where 
$   \maxuh \le \underline{t} \le \overline{t} \le -1/\minuh. $

\end{theorem}
\begin{proof} The proof follows the same arguments as the proof of Theorem \ref{th:threeway}.
\end{proof}

\begin{remark}
    Note that, when $N_2=N_2^\top \succeq 0$ has rank $2$, then the matrix $N_2$ is completely positive \cite{BermanP} and when choosing $\underline{t}=\overline{t}$, Theorem \ref{th:threewsym} then specializes to the result obtained in \cite{KG12}, using similar ideas.
\end{remark}
\begin{corollary} Theorem \ref{th:threeway} can also be reformulated using the arc-tangent angles 
corresponding to $\underline{t}_1$, $\overline{t}_1$, $\underline{t}_2$, $\overline{t}_2$, which yields the constraints~: 
\begin{eqnarray*}
\max_j \phi_{j} \le \underline{\alpha}_1 \le  \overline{\alpha}_1 \le \min_i \psi_{i}+\pi/2, \\ \max_i \psi_i \le \underline{\alpha}_2 \le  \overline{\alpha}_2 \le \min_j \phi_{j}+\pi/2.
\end{eqnarray*}
The factorization $N_2=LMR^\top$ then becomes
\begin{equation} \label{eq:12}
N_2=D_u\left[\cos(\overline{\alpha}_1\1-{\bf \Psi}), \sin(\underline{\alpha}_2\1-{\bf \Psi})\right]M
 \left[\cos(\overline{\alpha}_2 \1 -{\bf \Phi}),\sin(\underline{\alpha}_1\1-{\bf \Phi})\right]^\top D_v,
 \end{equation}
 where $M$ is nonnegative and is given by
 $$ M=\left[\begin{array}{cc}
   \cos(\overline{\alpha}_2-\overline{\alpha}_1)   & \sin(\underline{\alpha}_2-\overline{\alpha}_2) \\
   \sin(\underline{\alpha}_1-\overline{\alpha}_1)   & \cos(\underline{\alpha}_2-\underline{\alpha}_1)
 \end{array}\right]^{-1}$$
 $$= \left[\begin{array}{cc}
   \cos(\underline{\alpha}_2-\underline{\alpha}_1)   & \sin(\overline{\alpha}_2-\underline{\alpha}_2) \\
   \sin(\overline{\alpha}_1-\underline{\alpha}_1)   & \cos(\overline{\alpha}_2-\overline{\alpha}_1)
 \end{array}\right]/\left(\cos(\overline{\alpha}_2-\underline{\alpha}_1)\cos(\underline{\alpha}_2-\overline{\alpha}_1)\right). $$
 \end{corollary}
 \begin{proof}
 Equation \eqref{eq:12} is obtained by choosing the diagonal scaling matrices $D_L$ and $D_R$ in \eqref{eq:LR} such that
 $$ T_LD_L=\left[\begin{array}{cc}
   \cos(\overline{\alpha}_1)   & \sin(\underline{\alpha}_2) \\
   \sin(\overline{\alpha}_1)   &  -\cos(\underline{\alpha}_2)
 \end{array}\right], \quad T_RD_R=\left[\begin{array}{cc}
   \cos(\overline{\alpha}_2)   & \sin(\underline{\alpha}_1) \\
   \sin(\overline{\alpha}_2)   &  -\cos(\underline{\alpha}_1)
 \end{array}\right],
 $$
 and then multiplying the factors of the three terms of \eqref{eq:LR}
 and using the adjoint formula for the inversion of $M$. 
 \end{proof}

\section{Approximate rank-$2$ NMF}\label{sec:algorithm}

Let us now go back to the general singular value decomposition of $N$ as given in \eqref{svd2}, assuming $\rank ~N \geq 3$,
and use the rank-$2$ SVD approximation 
\[ M_2:=\hat \bu_1\hat\bv_1^\top + \hat\bu_2\hat\bv_2^\top\]
as the first rank-two approximation. Its error is given by
\begin{equation} \label{eq:bounds3}
 \|N-M_2\|_2 = \sigma_3, \quad  \|N-M_2\|_F = \sqrt{\sum_{i=3}^r \sigma_i^2}.
\end{equation}
We make the simplifying assumption that $M_2$ is unique, which is the case for the Frobenius norm when $\sigma_3>\sigma_2$. 
If not, the following discussion is still applicable to any arbitrary rank-$2$ approximation $M_2$. 
Recall also that, without loss of generality, we may assume that $N$ has no rows of zeros and that (thanks to Theorem \ref{thm:directsum}) $N$, and hence also $M_2$, has strictly positive dominant singular vectors.

If $M_2$ is nonnegative then this is obviously also the best nonnegative rank-two approximation and the bounds \eqref{eq:bounds3} are then met exactly. 
For example, the a-priori bound $\sigma_3 < \min_{i,j} N_{ij}$ is a sufficient condition for $M_2 \geq 0$, as proved in \cite[Lemma 6.7]{Gillis}.

If $M_2$ is not nonnegative, we can try to approximate $M_2$ by a nonnegative approximation $N_2$. Because of the triangle inequality of norms we then have
$$ \|N-N_2\| \le \|N-M_2\| + \|M_2-N_2\|,$$
which hopefully yields a good (but quite possibly suboptimal) nonnegative rank-two approximation $N_2$ to $N$. 
In Section \ref{sec:ANLS}, we propose to use $N_2$ as a starting
approximation for an iterative scheme, such as Alternating Nonnegative Least Squares (ANLS)\cite{Gillis}, to compute a hopefully better nonnegative rank 2 matrix approximation for the matrix $N$; in Section \ref{sec:numexp}, we present numerical experiments that successfully tested this approach. For now, we discuss below the computation of the solution for the modified problem $\min_{N_2\ge 0}\|N_2-M_2\|_F$. 
The benefit of considering this problem instead is that the matrix $M_2$ can now be written in the angular coordinates as in \eqref{eq:DcosD}.
Turns out, that a slightly relaxed version of the objective $\| N_2-M_2\|_F$ can be solved exactly in these coordinates.

\subsection{A relaxation for $\min\|N_2-M_2\|_F$}

We now turn to the geometric parametrization of $N_2$, as this will lead (after a relaxation) to a simple algorithm. In other words, let us use the factored forms from \eqref{eq:DPhiPsi}
$$ M_2=:L R^\top:=D_u [\cos({\bf \Psi}],\sin({\bf \Psi})] [\cos({\bf \Phi}],\sin({\bf \Phi})]^\top D_v$$
$$ N_2=:\widehat L \widehat R^\top :=\widehat D_u [\cos(\widehat {\bf \Psi}],\sin(\widehat {\bf \Psi})] [\cos(\widehat {\bf \Phi}],\sin(\widehat {\bf \Phi})]^\top \widehat D_v.
$$
By Lemma \eqref{lem:alpha}, the matrix $N_2$ is nonnegative if and only if there exist parameters $\theta_1,\theta_2 \in ]0,\pi/2[$ such that its angles $\psi_i$ and $\phi_j$ lie in the intervals
$$ \theta_{1} -\pi/2 \le \psi_i\le \theta_{2}, $$ 
$$ \theta_{2} -\pi/2 \le \phi_j \le \theta_{1}. $$ 
In terms of these coordinates, the optimization problem takes the form
\begin{align*}
    \min_{\theta_u,\theta_v,\widehat D_u,\widehat D_v,\widehat {\bf \Psi},\widehat {\bf \Phi}}~
&\|N_2-M_2\|_F \quad  \mathrm{where} 
\quad  N_2=\widehat D_u \cos(\widehat{\bf \Psi} 1_n^\top-\1_m\widehat {\bf \Phi}^\top)\widehat D_v \\
\mathrm{s.t.}~&(\theta_{1}-\pi/2) \1_m \le \widehat {\bf \Psi} \le \theta_{2}\1_m, \quad
(\theta_{2}-\pi/2) \1_n \le \widehat {\bf \Phi} \le \theta_{1}\1_n.
\end{align*}

\begin{remark}
    In the symmetric case, only one set of constraint is required~:
    $$  \theta -\pi/2\le \widehat\phi_i \le \theta,~\forall i \in [n]
    $$
    and in the cost function, only half of the variables are present. 
\end{remark}

Unfortunately, it does not appear obvious how to solve the nonlinear minimization problem above. Thus, we now show how to relax the nonlinear formulation to a more tractable problem. 
We will look for a cheaply computable, albeit possibly suboptimal, nonnegative solution for the minimization of $\|N_2-M_2\|_F$, with the idea to use it a starting point for a subsequent algorithm that aims at minimizing of $\|N_2-N\|_F$.

Namely, we consider the following bound on the objective function 
$$ \|N_2-M_2\|_F = \| LR^\top-\widehat L \widehat R^\top \|_F \le
\|(L-\widehat L)R^\top\|_F+\|\widehat L(R-\widehat R)^\top\|_F $$
$$\le \max(\|\widehat L\|_F,\|R\|_F) (\|L-\widehat L\|_F+\|R-\widehat R\|_F).
$$
As pointed out in Remark \ref{rem:balanced}, the coefficient
$\max(\|\widehat L\|_F,\|R\|_F)\approx \sqrt{\sigma_1+\sigma_2}$ cannot change too much, and therefore we may consider
\[ \min_{\widehat L \widehat R^\top \geq 0}\|L-\widehat L\|_F+\|R-\widehat R\|_F. \] 
Again, this is not easy to solve, and as a final relaxation (note that the norms get squared) we focus on 
\[ \min_{\widehat L \widehat R^\top \geq 0}\|L-\widehat L\|^2_F+\|R-\widehat R\|^2_F. \] 
The minimization problem over the variables $\widehat D_u$, $\widehat D_v$, $\widehat {\bf \Psi}$ and $\widehat {\bf \Phi}$ then
takes the form
$$
\|D_u [\cos({\bf \Psi}],\sin({\bf \Psi})]-\widehat D_u [\cos(\widehat {\bf \Psi}], \sin(\widehat {\bf \Psi})]\|^2_F $$
\begin{equation} \label{eq:min12}  +
\|D_v [\cos({\bf \Phi}],\sin({\bf \Phi})]-\widehat D_v [\cos(\widehat {\bf \Phi}], \sin(\widehat {\bf \Phi})]\|^2_F.
\end{equation} 
In \eqref{eq:min12}, angles and scalings can be tackled separately.

By solving a one-dimensional least squares problem, we get that for fixed $\widehat \Psi$, $\widehat \Phi$, the minimization problem over the scalings $\widehat D_u$ and $\widehat D_v$ can be explicitly solved term by term as
$$
\min_{\hat d_i^{(u)}}\|d_i^{(u)}[\cos(\psi_i),\sin(\psi_i)]-\hat d_i^{(u)}[\cos(\hat \psi_i),\sin(\hat \psi_i)]\|_2^2 = [d_i^{(u)}\sin(\psi_i-\hat\psi_i)]^2,
$$
$$
\min_{\hat d_j^{(v)}}\|d_j^{(v)}[\cos(\phi_j),\sin(\phi_j)]-\hat d_j^{(v)} [\cos(\hat \phi_j),\sin(\hat \phi_j)]\|_2^2 = [d_j^{(v)}\sin(\phi_j-\hat\phi_j)]^2.
$$
The corresponding optimal scale factors are
$$ \widehat D_u = D_u \cos(\Psi-\widehat \Psi) \quad \mathrm{and} \quad \widehat D_v = D_v \cos(\Phi-\widehat \Phi). $$
This allows us to omit the factors $\widehat D_u$ and $\widehat D_v$ from the problem.
In fact, we may omit the factors $\widehat \Psi$ and $\widehat \Phi$ as well, since the optimal choice is always determined directly by $\theta_1$ and $\theta_2$ as follows:
\begin{align} \label{eq:solhatpsi}
    \hat \psi_i = \begin{cases}
        \psi_i, &\text{if } \psi_i \in [\theta_1 - \pi/2, \theta_2], \\
        \theta_1 - \pi/2, &\text{if } \psi_i < \theta_1 - \pi/2, \\
        \theta_2,&\text{if } \psi_i > \theta_2,
    \end{cases}
\end{align}
and similarly for $\hat \phi_j$. 
From now on, we drop the parameters $\widehat D_u$, $\widehat D_v$, $\widehat \Psi$, and $\widehat \Phi$, and optimize over $\theta_1,\theta_2$.

Using \eqref{eq:solhatpsi}, we may restate the objective function. Each $\hat \psi_i$  (and analogously for $\hat \phi_j$) contributes to the cost by
\begin{align*}
\begin{cases}
    0, &\text{if } \psi_i \in [\theta_1 - \pi/2, \theta_2], \\
    [d_i^{(u)} \sin(\theta_2 - \pi/2 - \psi_i) ]^2, &\text{if } \psi_i < \theta_1 - \pi/2, \\
    [d_i^{(u)} \sin(\psi_i - \theta_1) ]^2, &\text{if } \psi_i > \theta_2.
\end{cases}
\end{align*}

Therefore, the total cost for \eqref{eq:min12} is
$f(\theta_1,\theta_2)=f_1(\theta_1)+f_2(\theta_2)$, where
\begin{equation}
 f_1(\theta_1)=\sum_{i= 1}^m \left( [d^{(u)}_i]^2 \sin^2(\max(0,\theta_{1}-\pi/2-\psi_i)) \right) +    \sum_{j=1}^n \left( [d^{(v)}_j]^2 \sin^2(\max(0,\phi_j-\theta_{1}) \right) \label{eq:anal1}
\end{equation}
and
\begin{equation*}
 f_2(\theta_2)=  \sum_{i=1}^m \left( [d^{(u)}_i]^2 \sin^2(\max(0,\psi_i-\theta_{2})) \right) + \sum_{j=1}^n \left([d^{(v)}_j]^2 [\sin^2(\max(0,\theta_{2}-\pi/2-\phi_j))\right). \end{equation*}

It is useful to note that we can in fact restrict the search to the subintervals
\begin{align}\label{eq:searchinterval}
    \min_i\psi_i +\pi/2 \le \theta_1 \le \max_j\phi_j, \quad
\min_j\phi_j +\pi/2 \le \theta_2 \le \max_i\psi_i,
\end{align} 
since it is easy to show that $f_1(\theta_1)$ and $f_2(\theta_2)$, respectively, monotonically increase outside these intervals with respect to inside them.

We aim to minimize $f(\theta_1,\theta_2)$ over $]0,\frac{\pi}{2}[^2$. For mathematical convenience, let us study $f$ over $[0,\frac{\pi}{2}]^2$; considering the closure of the domain mathematically simplifies our analysis and, as we shall later see in Lemma \ref{lem:unimodal}, is irrelevant for the purpose of computing its argument minimum.
In addition, we note the symmetry
$$f_2(\theta_2;{\bf \Psi},{\bf \Phi},D_u,D_v)=f_1(\theta_2;{\bf \Phi},{\bf \Psi},D_v,D_u),$$
implying that it suffices to provide a univariate optimization algorithm capable to compute $$\theta^*=\mathrm{argmin}_{\theta \in [0,\frac{\pi}{2}]} f_1(\theta;{\bf \Psi},{\bf \Phi},D_u,D_v),$$
for every possible choice of the input parameters ${\bf \Psi},{\bf \Phi},D_u,D_v$. For this reason, we now focus on this problem.

It is easy to verify that the function $\theta \in [0,\frac{\pi}{2}] \mapsto f_1(\theta)$, defined as in \eqref{eq:anal1}, is continuously differentiable on the closed interval\footnote{Here and throughout the paper, for a $k$-times differentiable function $g(\theta)$ defined on a closed interval $[a,b]$, we tacitly understand that its $j$-th derivative, $j \leq k$, at the endpoint $a$ (resp. $b$) is defined by using right (resp. left) one-sided limits. Hence, for example, $$f_1'(0)=\lim_{\varepsilon \rightarrow 0^+} \frac{f_1(\varepsilon)-f_1(0)}{\varepsilon}$$ is well defined in this sense.} $[0,\frac{\pi}{2}]$.
However, its second derivative is not continuous, which casts dark shadows over the applicability of standard higher-order optimization tools such as Newton's method. Instead, we provide an ad hoc algorithm to compute exactly a global minimizer of $f_1(\theta)$. 
Lemma \ref{lem:piecewiseformula} below is based on the idea of partitioning 
\begin{equation}\label{eq:partition}
    \left[0,\frac{\pi}{2}\right] = \bigcup_{k=0}^{q-1} [\theta_k,\theta_{k+1}], \qquad 0=\theta_0 < \theta_1 < \dots < \theta_q=\frac{\pi}{2},
\end{equation}
 where $\{\theta_k\}_{k=1}^{q-1} \subsetneq [0,\frac{\pi}{2}]$ denote the points of non-smoothness of $f_1$; in particular, $q \leq m+n+1$. Indeed, a point of non-smoothness is either an entry of the vector ${\bf \Phi}$ or an entry of the vector ${\bf \Psi} + \frac{\pi}{2} {\bf 1}$ (that is, an entry of ${\bf \Psi}$ shifted by $\frac{\pi}{2}$). However, depending on the angular parameters ${\bf \Phi}$ and ${\bf \Psi}$, it could happen that $q$ is (even significantly) smaller than its upper bound.  
 \begin{lemma} \label{lem:piecewiseformula}
   Let $f_1(\theta)$ be defined as in \eqref{eq:anal1}, and let $\{\theta_k\}_{k=0}^q$ define the partition \eqref{eq:partition}.
    For $k=0,\dots,q-1$, consider the functions
    \[ \theta \in [\theta_k,\theta_{k+1}] \mapsto g_k(\theta)= \sum_{\psi_i < \theta - \frac{\pi}{2}} [d^{(u)}_i \sin(\theta -\frac{\pi}{2} - \psi_i)]^2 + \sum_{\phi_j > \theta}[d^{(v)}_j \sin(\phi_j - \theta)]^2, \]
    and suppose further that, for all $k$, $g_k(\theta)$ is not identically zero in its domain of definition.
Furthermore, define
\begin{align}
        &a_k := \sum_{\psi_i < \theta_k - \pi/2} (d_i^{(u)})^2 (1 - 2 \sin^2 \psi_i) + \sum_{\phi_j > \theta_k } (d_j^{(v)})^2 (1 - 2 \cos^2 \phi_j), \notag \\ 
        &0 > b_k := \sum_{\psi_i <  \theta_k - \pi/2} (d_i^u)^2 \sin 2\psi_i -  \sum_{\phi_j > \theta_k} (d_j^v)^2 \sin 2\phi_j \label{eq:bi}
    \end{align}
 If a point $\theta_k^*$ exists in $[\theta_k,\theta_{k+1}] \subseteq [0,\frac{\pi}{2}]$ such that $g_k'(\theta_k^*)=0$, then
    \[ \theta^*_k =\begin{cases}
         \mathrm{sarctan}(\frac{b_k}{a_k})/2 &\mathrm{if } \ a_k\neq 0, \\
        \pi/4 &\mathrm{if} \ a_k=0,
    \end{cases} \]
    where 
    \[  \mathrm{sarctan}(x) := \begin{cases}
        \pi + \arctan(x) \ &\mathrm{if} \ x < 0,\\
        \arctan(x) \ &\mathrm{if} \ x > 0.
    \end{cases}      \]    
    Moreover, every global minimizer $\theta^*$ of $f_1(\theta)$ in $[0,\frac{\pi}{2}]$ does not lie at either endpoint, and satisfies $\theta^* \in \{\theta_k^*\}_{k=0}^{q-1}$.
\end{lemma}
\begin{proof}
Note first that $g_k(\theta)$ are all smooth, and in fact real-analytic, functions in their interval of definitions. Furthermore, by construction, $f_1(\theta)=g_k(\theta)$ on $[\theta_k,\theta_{k+1}]$.

We claim that $b_k < 0$, which proves \eqref{eq:bi} and justifies that the formulae defining $\theta_k^*$ are always well defined. Indeed, there  must be at least one term in the sum defining $b_k$, for otherwise $a_k=b_k=0$ and $g_k(\theta) \equiv 0$, contradicting the assumptions. Moreover,
    \[ -\pi/2 < \psi_i < \theta - \pi/2 \implies -\pi/2 < \psi_i < 0 \implies \sin 2\psi_i < 0, \]
    and similarly $\sin(2\phi_j) > 0$ for all $j \in [n]$.

    By exploiting basic trigonometric identities and some straightforward algebraic manipulations, we can show that the identities
    \[ f_1(\theta)=g_k(\theta) = a_k \cos^2(\theta) + \frac{b_k\sin(2\theta)}{2}, \quad f_1'(\theta)=g_k'(\theta) =  b_k \cos(2\theta)-a_k \sin(2\theta)  \]
    hold on the interval $[\theta_k,\theta_{k+1}]$. In particular, since $b_k<0$, $g_k$ is not a constant function on its interval of definition. In addition, these expressions immediately yield the stated formulae for $\theta_k^*$.

   Let now $\theta^*$ be an argument minimum of $f_1$ over $\theta \in [0,\frac{\pi}{2}]$; the existence of $\theta^*$ is clear because $f_1$ is continuously differentiable, and thus uniformly continuous, on $[0,\frac{\pi}{2}]$. As $f_1'(0)=b_1<0$ and $f'_1(\frac{\pi}{2})=-b_{q-1}>0$, implying that $\theta^* \not\in \{0, \frac{\pi}{2}\}$, it must happen that $\theta^*$ lies in the open interval $]0,\frac{\pi}{2}[$, and hence $f'_1(\theta^*)=0$.  
   On the other hand, necessarily $\theta^* \in [\theta_k,\theta_{k+1}]$ for some index $k$, implying in turn $g_k'(\theta^*)=0$ and hence $\theta^*=\theta_k^*$. 
\end{proof}
In the statement of Lemma \ref{lem:piecewiseformula}, the assumption that $g_k$ are not identically zero is not restrictive for the purposes of designing a numerical algorithm. Indeed, if $g_k \equiv 0$, then $f_1(\theta)=0$ for all $\theta \in [\theta_k,\theta_{k+1}]$, and therefore any element of the corresponding subinterval is a global minimizer.

One practical benefit of Lemma \ref{lem:piecewiseformula} is that the parameter $a_{k+1}$ can be computed in $\Theta(1)$ flops given the value of $a_k$:
\[ a_{k+1} := \begin{cases}
    a_k + (D_i^u)^2(1-2\sin^2 \psi_i) &\text{if } \theta_{k+1} = \psi_i+\pi/2 \\
    a_k - (D_j^v)^2(1-2\cos^2 \phi_j) &\text{if } \theta_{k+1} = \phi_j. \\
\end{cases}\]
An analogous observation holds for $b_k$, which implies that all $\theta_k^*$ can be computed  numerically in $\Theta(q)$ flops.
However, the angles must also first be ordered, which costs $\Theta(q\log(q))$ flops.
The bounds \eqref{eq:searchinterval} on the search interval further reduce the number of intervals $[\theta_k,\theta_{k+1}]$ to be considered, and thus the amount of angles to be sorted. In our numerical experiments (see Section \ref{sec:numexp}), this observation typically led to a very significant further speed-up. 

Finally, we note in Lemma \ref{lem:unimodal} that $f_1$ must be unimodal, implying that we can stop our search as soon as we have found one value $\theta_k^* \in [\theta_k,\theta_{k+1}]$. Once again, this can further improve the computational complexity.

\begin{lemma} \label{lem:unimodal}
    Under the same assumptions of Lemma \ref{lem:piecewiseformula}, the objective function $f_1$ as in \eqref{eq:anal1} is unimodal on $[0,\frac{\pi}{2}]$. In particular, $f_1$ has a unique stationary point $\theta^*\in ]0,\frac{\pi}{2}[$, which is also the global minimizer of $f_1$ on $[0,\frac{\pi}{2}]$.
\end{lemma}

\begin{proof}
We claim that $f_1'$ has at most one root in $[0,\frac{\pi}{2}]$. On the other hand, by the intermediate value theorem $f_1'$ has at least one root in $[0,\frac{\pi}{2}]$, because $f'_1(0)<0<f_1'(\frac{\pi}{2})$. (To check the latter inequalities, note that the constants $b_k$ are all negative by \eqref{eq:bi}.)  Thus, $f_1$ has precisely one stationary point $\theta^* \in ]0,\frac{\pi}{2}[$, and moreover $\theta^*$ must be a (local and global) minimizer of $f_1$.

    It remains to prove the claim. To this goal, using the same notation as in Lemma \ref{lem:piecewiseformula}, let us inspect the signs of the second derivatives $g_k''(\theta_k^*)$, assuming that $\theta_k^* \in [\theta_k,\theta_{k+1}]$. Explicitly, we have
    \begin{align*}
       g_k''(\theta) = -2a_k \cos(2\theta) - 2b_k \sin(2\theta).
       \end{align*}
    If $\theta_k^* = \pi/4$, the first addend vanishes and the second addend is positive. 
    Otherwise we may equivalently write
    \begin{align*}
        g_k''(\theta_k^*) &= \sin(2\theta_k^*) \left( \frac{-2a_k}{\tan(2\theta_k^*)} - 2b_k \right) = \sin(2\theta_k^*) \left( -\frac{2a_k^2}{b_k} - 2b_k \right) \\
        &= \frac{\sin(2\theta_k^*)}{-b_i} (2a_k^2 + 2b_k^2).
    \end{align*}  
    Since $\theta_k^* \in ]0,\pi/2[$, we conclude that $g_k''(\theta_k^*)>0$.

    By Lemma \ref{lem:piecewiseformula}, any root of $f_1'$ in $[0,\frac{\pi}{2}]$ must be equal to $\theta_k^* \in ]0,\frac{\pi}{2}[$ for some index $k$ such that $\theta_k^* \in [\theta_k,\theta_{k+1}]$. There are two possible cases to consider:
    \begin{enumerate}
        \item If $\theta_k^* \in ]\theta_k,\theta_{k+1}[$, then $f_1''(\theta_k^*)=g_i''(\theta_k^*) > 0$.
        \item If instead $\theta_k^* \in \{\theta_k,\theta_{k+1}\}$, then $f_1$ is not twice differentiable there; however, in this case, the left and right second derivatives of $f_1$ at $\theta_k^*$ coincide with the second derivatives at $\theta_j^*$ of two consecutive $g_j(\theta),g_{j+1}(\theta)$, which are both positive by the same argument.
    \end{enumerate} 
 We conclude that (in either case and for all $k=0,\dots,q-1$) the a.e. differentiable function $f_1'$ is strictly increasing in every small enough open ball centered at $\theta_k^*$. Therefore, $f_1'$ can have at most one root in $[0,\frac{\pi}{2}]$. This concludes the proof.
    \end{proof}

The procedure is described in Algorithm \ref{alg:quadrant}. 
Its computational cost is dominated by the rank-2 SVD at the beginning, which can be computed in $\Theta(mn)$ flops.

\begin{algorithm}[H]
    \KwData{Nonnegative matrix $N \in \RR_+^{m\times n}$.}
    \KwResult{Nonnegative factors $L_+ \in \RR_+^{m\times 2}$ and $R_+ \in \RR_+^{n \times 2}$ such that $N \approx L_+R_+^\top$.}
    
    $\widehat U, \widehat V \gets \mathrm{svd}(N,2)$; 
    Remove possible zero rows from $\widehat U$ and $\widehat V$.

    $D_u \gets \mathrm{diagm}( \|\widehat U(i,:)\|_2^{-1})$; $D_v \gets \mathrm{diagm}( \|\widehat V(i,:)\|_2^{-1})$;

    ${\bf \Psi} \gets \arctan(\widehat U(:,1),\widehat U(:,2))$;   ${\bf \Phi} \gets \arctan(\widehat V(:,1),\widehat V(:,2))$;

    Solve the optimal angles $\theta_1,\theta_2$ using Lemma \ref{lem:piecewiseformula}. 

    $\widehat {\bf \Psi} \gets \min(\max({\bf \Psi},\theta_1-\frac{\pi}{2}),\theta_2)$;  $\widehat {\bf \Phi} \gets \min(\max({\bf \Phi},\theta_2-\frac{\pi}{2}),\theta_1)$;

    $\widehat D_u \gets D_u \cos^2({\bf \Psi} - \widehat {\bf \Psi})$;   $\widehat D_v \gets D_v \cos^2({\bf \Phi} - \widehat {\bf \Phi})$;   

    Find the exact NMF-2 for $\widehat D_u \cos(\widehat {\bf \Psi} - \widehat {\bf \Phi})\widehat  D_v =: L_+R_+^\top$ using \eqref{eq:factors}.
    
    \caption{Quadrant (QDR)} \label{alg:quadrant}
\end{algorithm}

\section{A new algorithm: SVD-based initialization for ANLS}\label{sec:ANLS}

The solution $N_2=L_+R_+^\top$ of the relaxed problem, that can be cheaply computed by Algorithm 1 described in Section \ref{sec:algorithm}, is not always competitive. However, it can be used as a starting value for finding the best nonnegative rank-two approximation of $N$, using a variant of the Alternating  Nonnegative Least Squares (ANLS) method \cite{Gillis}. As we will illustrate in the numerical experiments in Section \ref{sec:numexp}, this modified starting point eventually provides several advantages with respect to the standard initialization of ANLS.

\subsection{The ANLS update and its symmetric version}

Assume we have a nonnegative rank-two approximation in factored form $N_2=L_+ R_+^\top$ of a given nonnegative matrix $N$. The ANLS method minimizes the Frobenius norm of the distance $N-L R^\top$
by alternating between updating $L$ and $R$, which can be decoupled
in subproblems where only one row of $L$ or $R$ is considered at a time \cite{Gillis}. For instance, when $L$ is fixed and $R$ is updated, the resulting constrained linear least squares problem is
\begin{equation*}
    \min_{{\bf x}_+\in \mathbb{R}_+^2} \| L {\bf x}_+  - {\bf b} \|_2, 
\end{equation*} 
where ${\bf b}\ge 0$ is the corresponding column of $N$. 
The complexity of one ``sweep" (that is, one update of $R$ followed by one update of $L$) is $\Theta(mn)$. At the same time, this setting clearly begs the question of how to initialize ANLS before the first sweep, where no previous output exists. Indeed, ANLS can be sensitive to the choice of its starting point \cite{Langville}.

We now briefly discuss how to possibly address the symmetric case, as this seems to not have received much attention in the literature. We start from a rank-$2$ nonnegative ``three way factorization" 
$N_2:=LML^\top$ and update the factor $L$ in order to minimize the Frobenius norm of the error $\|N-LML^\top\|_F^2$. The proposed scheme updates each column of the right factor $L^T$ and then also updates the left factor such that the factorization remains symmetric. 
Let us consider updating only the last column $x$ of the matrix $L^\top$, where we partition the matrices as follows
$$ LML^\top =: 
\left[\begin{array}{cc}  L_1 \\ {\bf x}^\top \end{array}\right]M
\left[\begin{array}{cc}  L_1^\top & {\bf x} \end{array}\right] , \qquad N=:\left[\begin{array}{cc}  N_{11} & {\bf b} \\ {\bf b}^\top & n \end{array}\right].
$$
Now solve exactly the one sided minimization problem
$$ \min_{\hat {\bf x}\in \R_+^2} \left\|
\left[\begin{array}{cc}  L_1 \\ {\bf x}^\top \end{array}\right]M
\left[\begin{array}{cc}  L_1^\top & \hat {\bf x} \end{array}\right] -
\left[\begin{array}{cc}  N_{11} & {\bf b} \\ {\bf b}^\top & n \end{array}\right]\right\|_F^2
$$
which has the same solution as the standard problem 
$$ \min_{\hat {\bf x}\in \R_+^2} \left\|
\left[\begin{array}{cc}  L_1 \\ {\bf x}^\top \end{array}\right]M
 \hat {\bf x}  - 
\left[\begin{array}{cc} {\bf b} \\ n \end{array}\right]\right\|_F^2
$$
and then replace both ${\bf x}$ and $\hat {\bf x}$ by their average $\frac12({\bf x}+\hat {\bf x})$. We now show that the symmetrized error and 
the unsymmetric error are of the same order.
Define
$$ h_1:=  2\frac{(\bl-\br)}{2}^\top \frac{(\bl-\br)}{2} =
(\bl^\top \bl +\br^\top \br) 
 -2\frac{(\bl+\br)}{2}^\top \frac{(\bl+\br)}{2} \geq 0,
$$
$$ h_2:= - \frac{({\bf x}-\hat {\bf x})}{2}^\top M\frac{({\bf x}-\hat {\bf x})}{2} = {\bf x}^\top M \hat {\bf x}  -\frac{({\bf x}+\hat {\bf x})}{2}^\top M\frac{({\bf x}+\hat {\bf x})}{2}.
$$
Specializing to $\bl:=L_1M{\bf x}-{\bf b}$ and $\br:=L_1M\hat {\bf x}-{\bf b}$, and denoting by $h(\cdot)$ the squared Frobenius norm of a vector, we then have
$$  h\left(\frac{{\bf x}+\hat {\bf x}}{2}\right) = h(\hat {\bf x}) -h_1 + h_2.
$$
Both $h_1$ and $h_2$ are higher order corrections when all errors are small, showing $h\left(\frac{{\bf x}+\hat {\bf x}}{2}\right) \approx h(\hat {\bf x})$.

An alternative approach is to solve the symmetric minimization problem successively for each element of ${\bf x}$. The corresponding cost function is quartic and its derivative is cubic, which implies 
one can solve for each element $x_i, \; i=1,2$ by
comparing the different solutions \cite{VanGZD16}.

\bigskip

\section{Numerical experiments}\label{sec:numexp}

In this section, we present several numerical tests whose objective is to compute, by means of the ANLS method, the rank-$2$ approximate NMF of a nonnegative matrix input $N$ whose rank-$2$ truncated SVD $M_2$ is not nonnegative (because otherwise the problem is trivial). We compare Algorithm \ref{alg:quadrant} with four other possible initialization methods for ANLS:
\begin{enumerate}
    \item The successive projection algorithm (SPA) first presented in \cite{AraEtAl01} and successfully used in \cite{GilKP14} for NMF-2 based hierarchical clustering. While SPA is provably optimal for the special subclass of nonnegative matrices studied in \cite{GilKP14}, this is not true for general input. 
For a fair comparison, we use SPA to initialize ANLS even though in \cite{GilKP14} it is used as the sole algorithm to compute a NMF.

    \item Nonnegative Double Singular Value Decomposition (NNDSVD), introduced in \cite{BG08}.

    \item Nonnegative SVD with low-rank correction (NNSVDLRC), developed in \cite{Atifetal19} and applied in \cite{Atifetal21}.
    The method works for general $r$, and for $r=2$ it takes a similar structure as NNDSVD.

    \item Randomly generated starting point.
    
\end{enumerate}

In the experiments, we used the implementation of nonnegative least squares given by Gillis, Kuang, and Park in \cite{GilKP14} as a subroutine for the ANLS algorithm. As a convergence test within the ANLS iterations, we measure the ratio
\begin{equation}\label{eq:residual}
     \mathrm{res} := \sum_{i=1}^m \frac{\| R^{(k+1)}_i - R^{(k)}_i\|_2}{\|R_i^{(k+1)}\|_2} + \sum_{j=1}^n\frac{\| L^{(k+1)}_j - L^{(k)}_j\|_2}{\|L_j^{(k+1)}\|_2},
\end{equation}
where $R^{(k)}_i$ denotes the $i$th row of the factor $R$ at the $k$th iteration, and similarly for the factor $L$.
We fix a tolerance $\epsilon$, and terminate the algorithm once either the residual in \eqref{eq:residual} is lower than $\epsilon$ or a maximum number of iterations is reached.

Throughout this section we denote the output of ANLS with the five different initialization methods described above by, respectively,  $N_{\mathrm{QDR}}$, $N_{\mathrm{SPA}}$, $N_{\mathrm{NNSVDLRC}}$, $N_{\mathrm{NNDSVD}}$, and $N_{\mathrm{rand}}$. 
We refer to the optimal rank-$2$ nonnegative approximation by
\[ N^* := \argmin_{\mathrm{rank}_+(X) \leq 2} \| N - X \|_F. \]
When the exact minimizer $N^*$ is not known, we refer to the output closest to the original matrix $N$ by $$N_{\min} \in \{ N_{\mathrm{QDR}}, N_{\mathrm{SPA}}, N_{\mathrm{NNSVDLRC}}, N_{\mathrm{NNDSVD}},  N_{\mathrm{rand}} \}.$$
In these cases, the relative performance of the approximation $N_{\mathrm{alg}}$ computed using the initialization method $\mathrm{alg}$ is measured by
\[ \delta_{\mathrm{alg}} := \frac{\| N - N_{\mathrm{alg}}\|_F}{\|N - N_{\min}\|_F} \geq 1.  \]

We were not able to find general-purpose benchmark tests for rank-$2$ NMF in the literature. We found some papers that do provide data, but they are primarily specialized on particular subclasses of nonnegative matrices, whereas we wish to test our general-purpose method on data that is not structured (beyond being nonnegative and such that their nearest rank-$2$ matrix is not already nonnegative -- as otherwise the problem is trivial). 
Thus, our test consist of synthetic data, i.e., randomly generated nonnegative matrices.
Specifically, we consider three classes of input:

\begin{enumerate}
    \item[A)] Dense full rank nonnegative $n\times n$ (A1) and $n \times 150$ (A2) matrices.
    \item[B)] Nonnegative matrices of rank $2$ with added noise.
    \item[C)] Nonnegative $4\times 4$ integer matrices for which we know the optimal $N^*$.
\end{enumerate}

The code for the experiments can be found in \url{https://github.com/EtnaLindy/NMF2_QDR}.
A summary of the outcomes is given in Table \ref{tab:summary}.

\begin{table}[h]
    \centering
    \scriptsize
    \begin{tabular}{|c|ccccc|}
    \hline
         Case A1. $n = 800$ & SPA & NNSVDLRC & NNDSVD & rand & QDR \\
mean time & 0.026 & 0.087 & 0.011 & 0.04 & \textbf{0.01}\\
max time & 0.438 & 0.503 & 0.386 & \textbf{0.381} & 0.392\\
mean acc & 1 & 1 & 1 & 1 & 1\\
max acc & \textbf{1} & 1.147 & \textbf{1} & 1.033 & \textbf{1}\\
mean acc init & \textbf{1} & 1.044 & \textbf{1} & 1.338 & \textbf{1}\\
max acc init & \textbf{1} & 1.622 & 1.033 & 10.337 & \textbf{1}\\ \hline\hline
Case A2. $n = 15000$ & SPA & NNSVDLRC & NNDSVD & rand & QDR \\
mean time & 0.022 & 0.133 & 0.051 & 0.241 & \textbf{0.013}\\
max time & 1.725 & 1.869 & 1.786 & 1.766 & \textbf{0.22}\\
mean acc & 1 & 1 & 1 & 1 & 1\\
max acc & \textbf{1} & 1.256 & \textbf{1} & \textbf{1} & \textbf{1}\\
mean acc init & \textbf{1} & 1.056 & \textbf{1} & 1.499 & \textbf{1}\\
max acc init & \textbf{1} & 1.577 & 1.003 & 17.836 & \textbf{1}\\ \hline\hline
Case B. $n = 100 $& SPA & NNSVDLRC & NNDSVD & rand & QDR \\
mean time & 0.004 & 0.006 & 0.005 & 0.003 & \textbf{0.001}\\
max time & 0.017 & 0.019 & 0.019 & 0.015 & \textbf{0.002}\\
mean acc & 1 & 1 & 1 & 1 & 1\\
max acc & 1 & 1 & 1 & 1 & 1\\
mean acc init & 1.008 & 1.009 & 1.162 & 3.565 & \textbf{1}\\
max acc init & 1.232 & 1.02 & 1.263 & 4.171 & \textbf{1.002}\\ \hline \hline
Case C. $n = 4$ & SPA & NNSVDLRC & NNDSVD & rand & QDR \\
mean time & 0.00088 & 0.00087 & 0.00099 & 0.00076 & \textbf{0.00049}\\
max time & 0.027 & 0.019 & 0.022 & 0.02 & \textbf{0.018}\\
mean acc & \textbf{1} & 1.001 & \textbf{1} & 1.058 & \textbf{1}\\
max acc & 1.145 & 2.051 & 1.083 & 12.609 & \textbf{1.004}\\
mean acc init & 1.12 & 1.231 & 1.721 & 5.357 & \textbf{1.02} \\
max acc init & \textbf{5.254} & 17.212 & 15.228 & 58.648 & 6.865\\ \hline
    \end{tabular}
    \caption{Summary of the numerical experiments described in Section \ref{sec:numexp}.}
    \label{tab:summary}
\end{table}

\subsection{Dense high-rank matrices} \label{sec:sampling}

Recall that we aim to generate a nonnegative matrix $N$ whose nearest rank-$2$ matrix $M_2$ is (with reasonably high probabiility) not nonnegative. 
To justify our choice of a probability distribution to generate our data, consider an $m \times n$ ($m \geq n$) random matrix $X$, whose entries are i.i.d. with support on $[0,+\infty[$, finite mean $\mu$, and finite variance $\eta^2$. 
Clearly, $X=\mu {\bf 1}_n {\bf 1}_m^{\top} + Y$ is the sum of a deterministic rank-$1$ matrix and a centered random matrix with the same variance $\eta^2$. 
According to universality results in random matrix theory \cite{chafai,TVK}, in the limit $m,n \rightarrow \infty$ and provided that $\eta \lesssim (\mu/2)\sqrt{n}$, we expect the largest two singular values of $X$ to behave as $\sigma_1(X) \approx \mu \sqrt{mn}$ and $\sigma_2(X) \approx \eta (\sqrt{m} + \sqrt{n})$. 
More informally, unless the standard deviation $\eta$ grows significantly faster than the mean $\mu$, the expectation $\mathbb{E}[X]=\mu {\bf 1}_n {\bf 1}_m^\top$ dominates the noise $Y$, creating a large singular value gap $\sigma_1/\sigma_2$. 
In practice, for many common distributions, these asymptotic predictions are likely to be accurate already for moderate values of $m,n$. 
Assuming $m = \alpha n$ for some constant $\alpha \geq 1$, we obtain the estimate
\[ \frac{\sigma_1}{\sigma_2} \approx \frac{\mu n \sqrt{\alpha}}{\eta(\sqrt{\alpha n} + \sqrt{n})}  = \Theta \left( \frac{\mu}{\eta} \sqrt{n} \right).  \] 
By Lemma \ref{lem:conditions}, we expect $M_2 \geq 0$ with high probability when $\sigma_1/\sigma_2$ is large. 
Many popular nonnegative distributions such as the $[0,1]$-uniform distribution (realized in MATLAB by the command \texttt{rand(m,n)}) or the folded normal distribution (\texttt{c*abs(randn(m,n))} in MATLAB, where $c>0$ is allowed to depend on $m$ and $n$) satisfy $\eta/\mu = \Theta(1) \Rightarrow \sigma_1/\sigma_2 = \Theta(\sqrt{n})$. 
As a result, they are unsuitable for our purposes unless $m,n$ are very small.  
To circumvent this issue, we seek instead a heavy-tailed probability distribution such that $\eta/\mu = \Theta(\sqrt{n}) \Rightarrow \sigma_1/\sigma_2 = \Theta(1)$. 
Our choice is the lognormal distribution defined by $e^Z$ where $Z \sim \mathcal{N}(0,\log m)$ is normal with mean $0$ and variance $\log m$. 
This is generated in MATLAB by \texttt{exp(sqrt(log(m))*randn(m,n))}. 
Still under the assumption $m=\alpha n$, for this distribution $\mu=\Theta(\sqrt{n})$ and $\eta=\Theta(n)$, thus achieving a small singular value gap $\sigma_1/\sigma_2$, as desired.

\begin{figure}[h!]
    \centering
\begin{subfigure}[t]{0.55\textwidth}
    \centering
    \includegraphics[width=\linewidth]{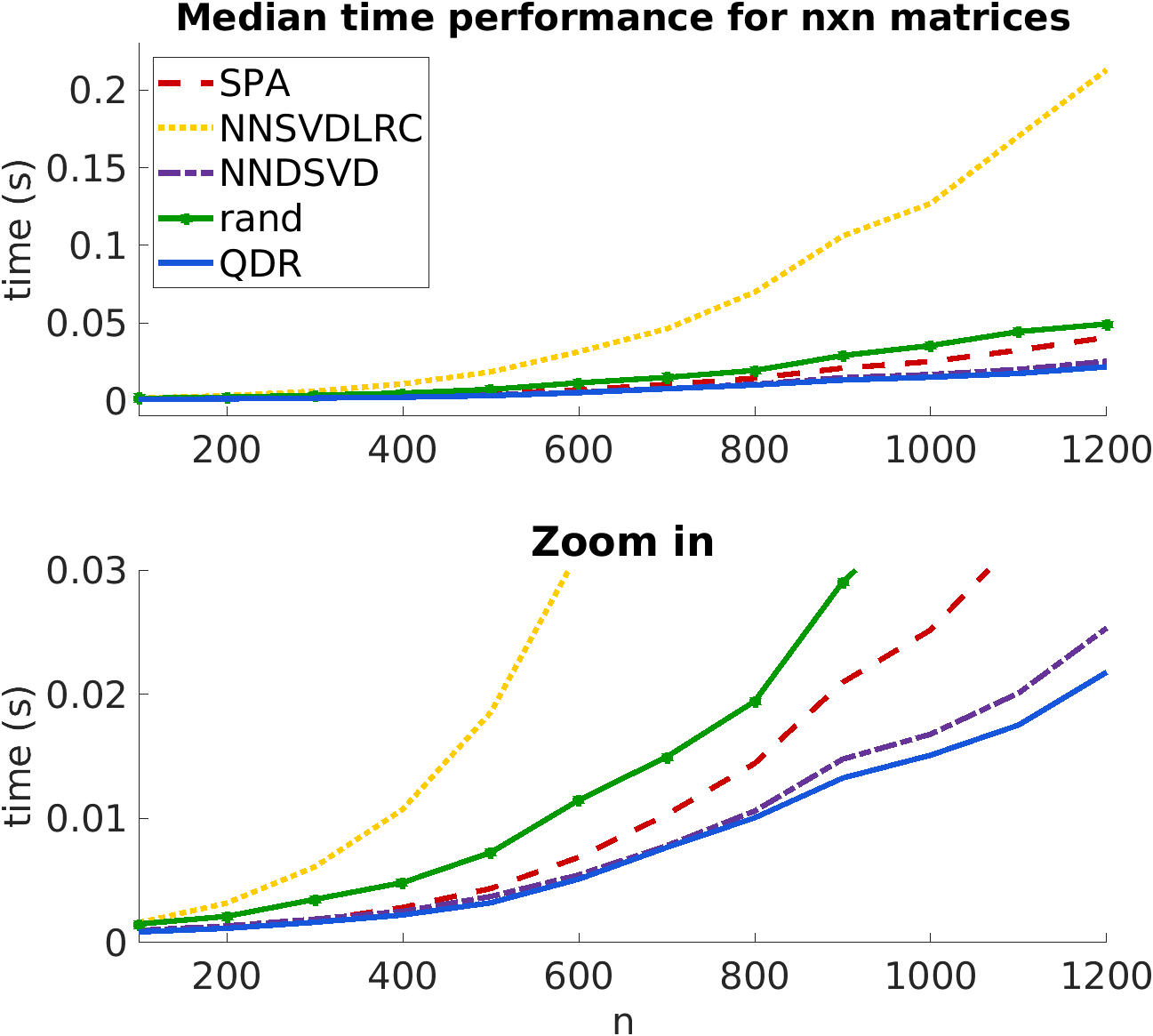}
\end{subfigure}
~
\begin{subfigure}[t]{0.42\textwidth}
    \centering
    \includegraphics[width=\linewidth]{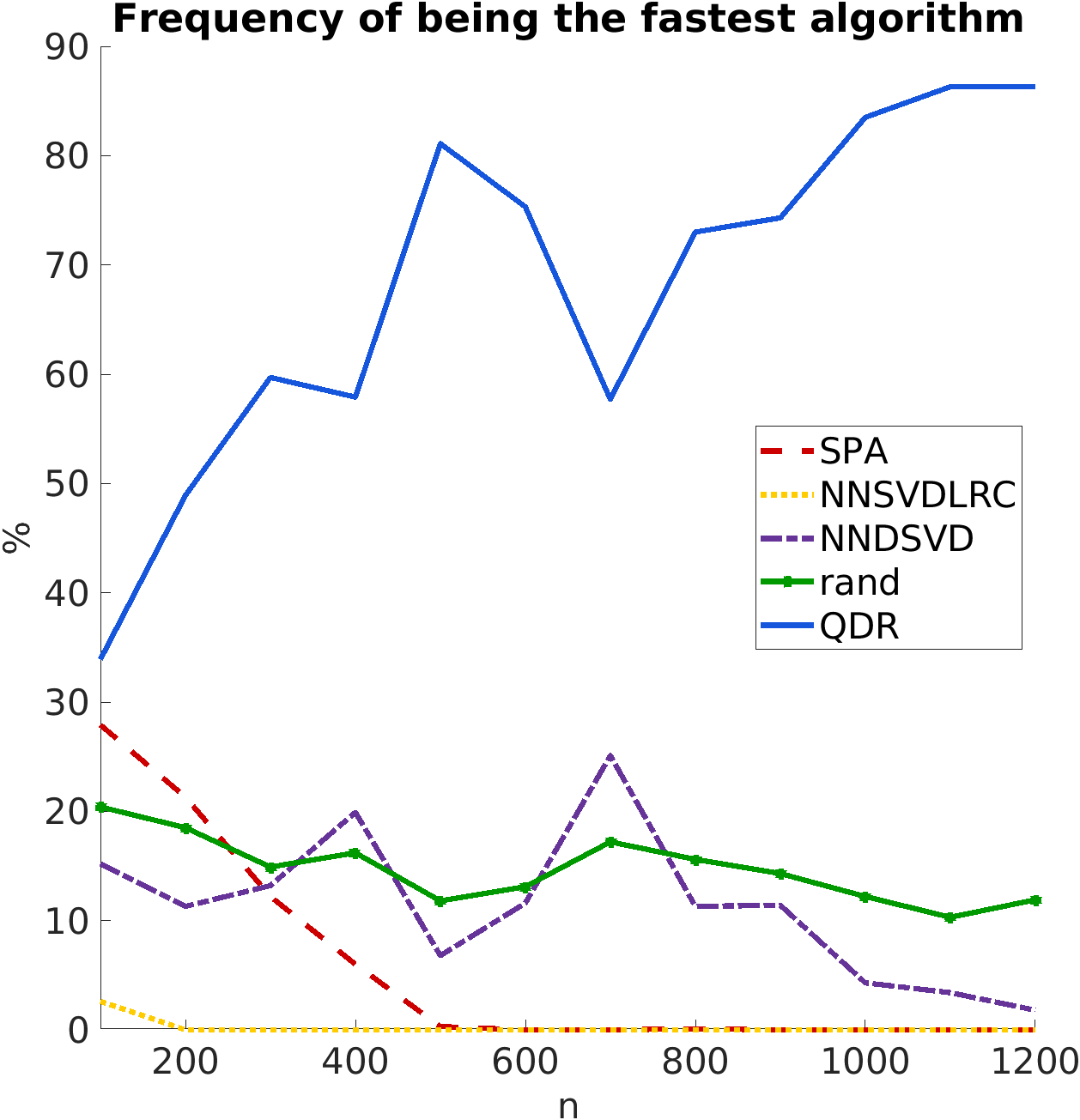}
\end{subfigure}
\caption{ANLS with convergence tolerance $\epsilon = 10^{-3}$ and maximum number of iterations $1000$ applied to $n\times n$ matrices sampled as described in Section \ref{sec:sampling}. On the right, one can see the percentage of cases where a specific initialization algorithm performed the fastest.}
\label{fig:lognormal_time}
\end{figure}

\begin{figure}[h!]
    \centering
    \begin{subfigure}[t]{0.48\textwidth}
    \centering
    \includegraphics[width=\linewidth]{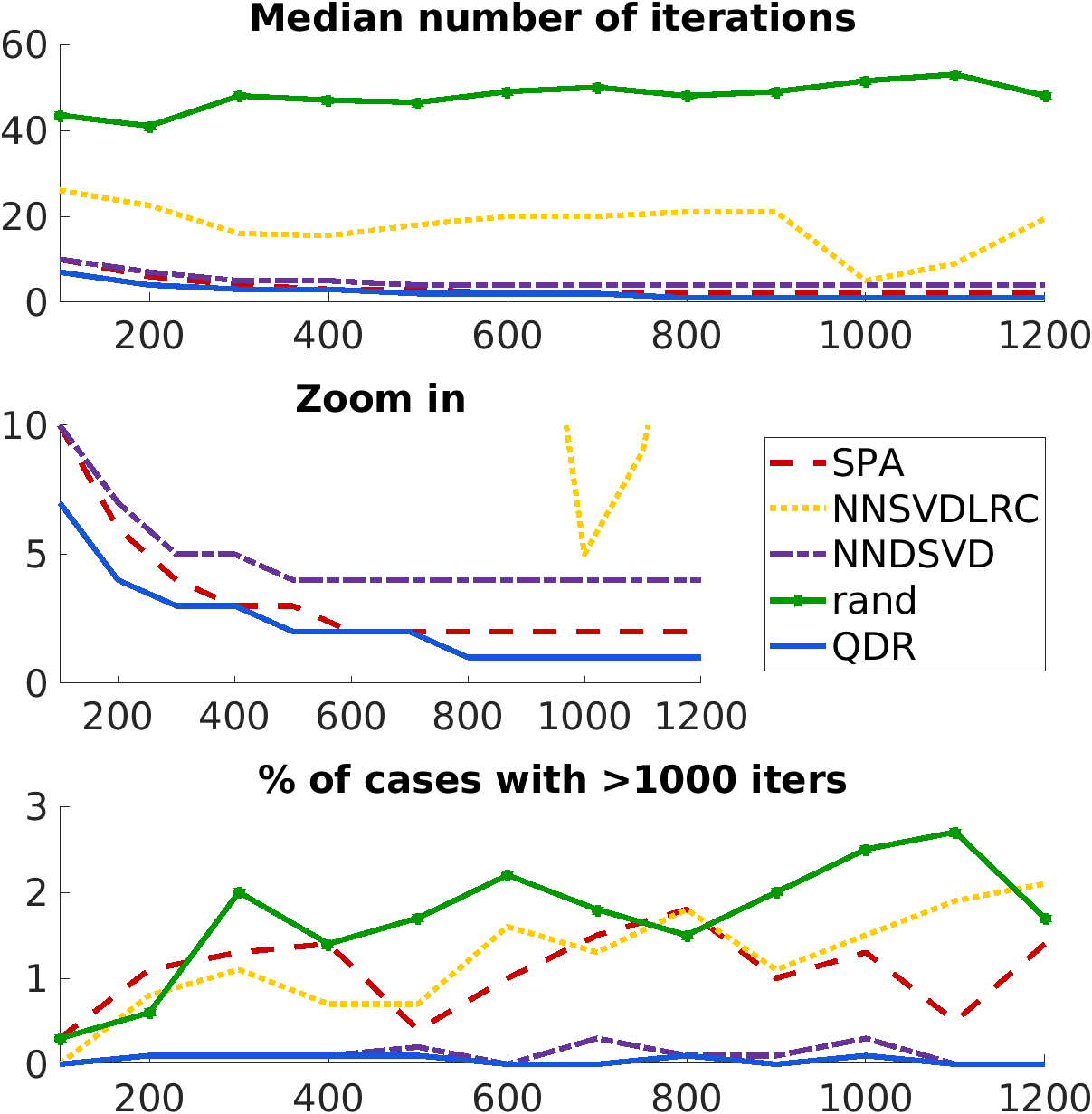}

\end{subfigure}~
    \begin{subfigure}[t]{0.48\textwidth}
    \centering
    \includegraphics[width=\linewidth]{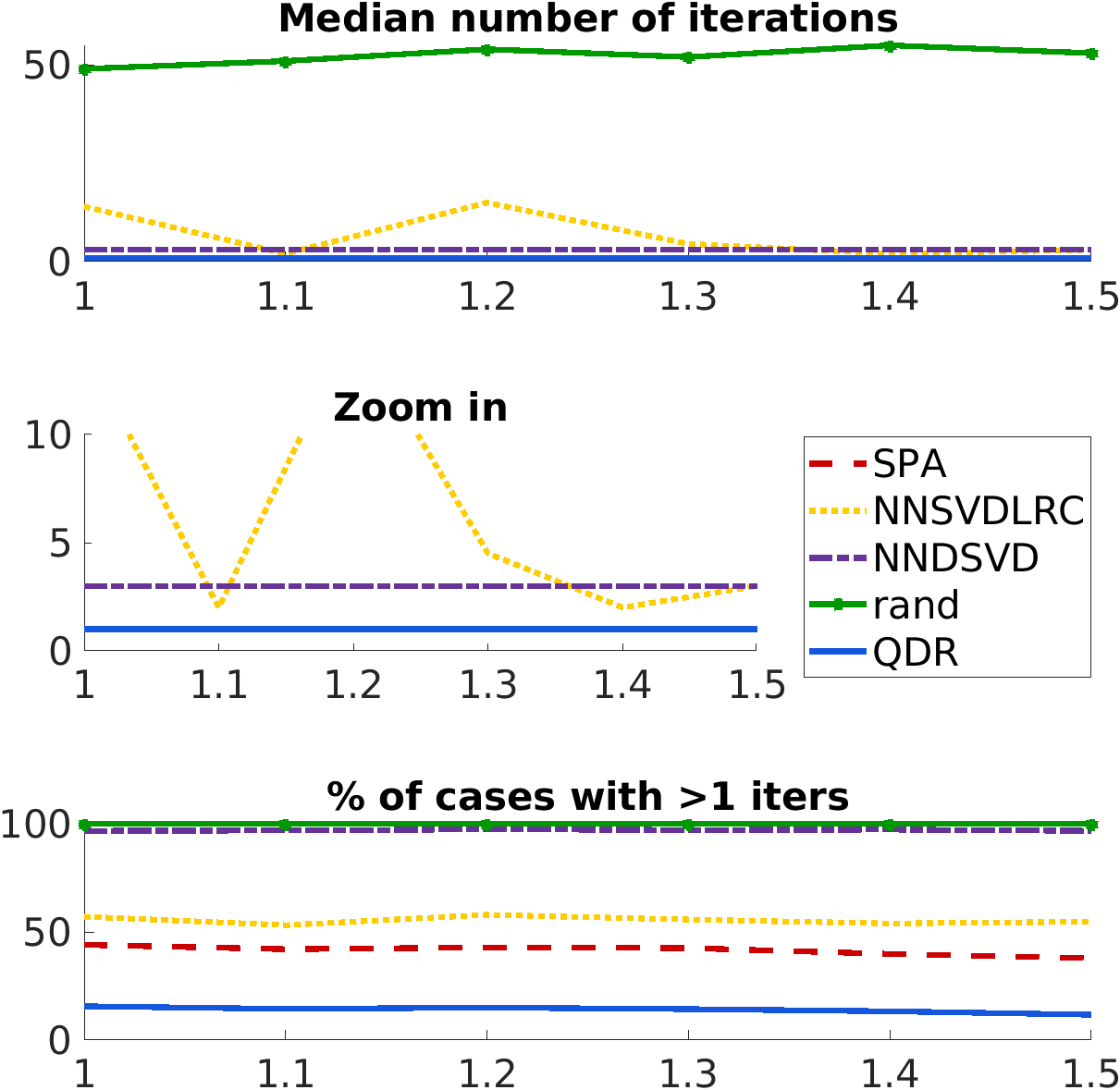}
\end{subfigure}
\caption{ANLS with convergence tolerance $\epsilon = 10^{-3}$ and maximum number of iterations $1000$ applied to $n\times n$ (left) and $n\times 150$ (right) matrices sampled as described in Section \ref{sec:sampling}. For $n\times 150$ matrices, initialization with QDR and SPA, the median number of ANLS iterations was equal to 1 for all $n$, which is why the frequency of encountering more than 1 iterations is counted.}
\label{fig:lognormal_iters}
\end{figure}

In Figure \ref{fig:lognormal_time}, we plot the performance of ANLS initialized with various algorithms, as the dimension varies. 
The QDR initialization achieved the best outcome, both in terms of the median running time and the number of cases where ANLS with QDR initialization was the fastest among the methods.
Figure \ref{fig:lognormal_iters} reports statistics on the number of iterations for randomly sampled square matrices (left) and tall matrices of size $n\times 150$ (right, $10^5 \leq n \leq 1.5 \cdot 10^5$).
In both experiments the convergence parameter was $\epsilon = 1e-3$ and the maximum number of iterations was $1000$.
The latter parameter is intentionally large, to ensure that in most cases the convergence criterion is met before reaching the maximum number of iterations. 

For the initializations with QDR, SPA, and NNDSVD, we observed a somewhat surprising decrease in the median number of ANLS iterations needed for convergence as the dimension grows. For matrices of size $n\times 150$, for all $n$ the median number of iterations for QDR and SPA was $1$, i.e., ANLS terminated immediately. 
In both tests, the number of ANLS iterations needed to reach the convergence criterion was consistently lowest when using QDR as the initialization method (number of cases with more than 1000 iterations and the median number of iterations).

\subsection{Rank-$2$ data with noise}

In this experiment, we wish to simulate a situation where the input $N$ is supposed to be rank $2$, but data have been corrupted by noise, and one aims to correct for this by performing an approximate rank-$2$ NMF. 

\begin{figure}[h!]
    \centering
    \includegraphics[width=\linewidth]{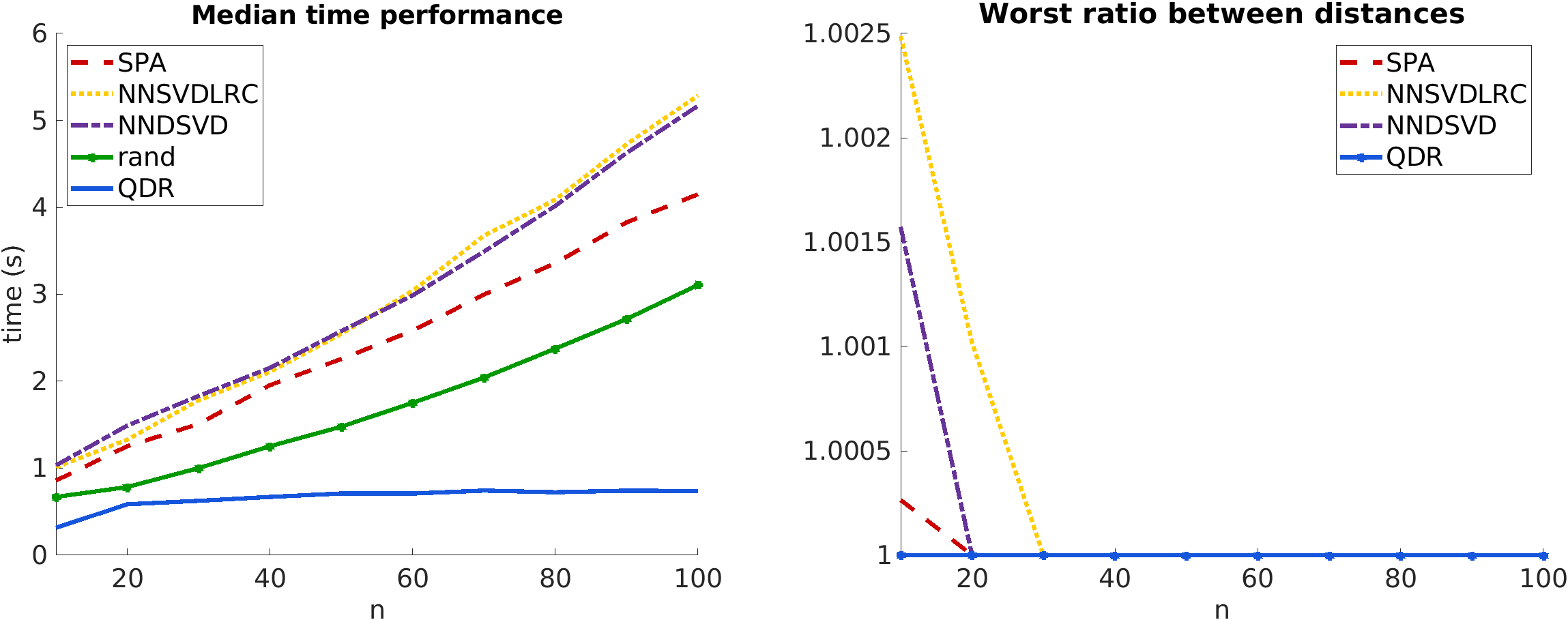}
    \caption{ANLS with maximum of 1000 iterations and $\epsilon = 10^{-3}$ applied to low-rank data sampled from the boundary with additive noise.  Random-point initialization has been omitted from the figure on the right.}
    \label{fig:lowdimnxn}
\end{figure}

We generate our input starting from \eqref{eq:DcosD}, setting $\phi_1=\dots=\phi_{ n/2}=0$, $\psi_m=\psi_{m-1}=\dots=\psi_{ m/2}=\pi/2$, and drawing rest of the angles uniformly from $[0,\frac{\pi}{2}]$. 
The weights $D_u$ and $D_v$ are uniformly distributed in $[0,1]$. 
This constructs a rank-$2$ nonnegative matrix on the boundary of the set of nonnegative rank $2$ matrices. 
Finally, we add noise distributed as a folded normal, keeping only outcomes whose rank-$2$ best approximation is not nonnegative. 
The reason to start the construction from a point on the boundary is to make sure that this happens with reasonably high probability. 
In Figure \ref{fig:lowdimnxn}, we display the performance of ANLS with various initialization methods. 
The experimental outcome favors our method; most of the time all the methods seem to converge to an equally good approximation, but initialization with Algorithm \ref{alg:quadrant} converges in the least amount of time.

\subsection{Small-size input with known exact solution}

For this experiment, we used a database of $4\times 4$ nonnegative integer matrices, constructed by Lindy in \cite{Lindy}. 
There, the elements of $10^5$ matrices were drawn uniformly from the set of integer nonnegative matrices whose entries sum up to $1000$, and then the trivial cases (those whose best rank-$2$ approximation is already nonnegative) were dismissed, yielding a total of 25357 ``interesting" test matrices. 
In \cite{Lindy}, all the critical points of the cost function $\|N-N_2\|_F$ were computed  using a homotopy continuation method together with algebraic-geometric techniques for finding the number of possible critical points for generic data $N \in \mathbb{N}^{4\times 4}$. 
Therefore, the database also provides the certified best approximations solutions $N^*$. 
We note that, while the approach in \cite{Lindy} guarantees to compute the optimal solution, it is computationally very demanding and hence impractical for larger input sizes.

In Table \ref{tab:pointofconvergence} one can see, for these test examples, how close to the true optimum $N^*$ the output of ANLS with starting point from Algorithm \ref{alg:quadrant} is. The same statistics have been reported for the other initialization methods as well. 
For all the starting points, ANLS was then performed (with $\epsilon = 10^{-5}$ and maximum of $10^4$ iterations).

In only 11 cases, ANLS with starting point from Algorithm \ref{alg:quadrant} converges to a point significantly far away from the exact optimum $N^*$ (Table \ref{tab:pointofconvergence}). In these 11 cases, the computed $N_{\mathrm{QDR}}$ was nevertheless a relatively good approximation, as the ratio
\[  \frac{\|N_{\mathrm{QDR}} - N\|_F}{\|N^* - N\|_F} \lesssim 1.004 \]
for all the sampled matrices.
In comparison, for the SPA starting point, the matrix $N_{\mathrm{SPA}}$ significantly differs from $N^*$ in 79 cases and the maximal ratio of distances was $1.14$, and for the NNDSVD starting point there were 114 cases with maximum ratio of distances 1.08.
For NNSVDLRC and for the random-point initialization, the number of cases that did not converge to the optimum were 258 and 1250, and the maximum ratios were 2.0514 and 12.6088 respectively.

\begin{figure}[h]
    \centering
    \begin{subfigure}[t]{0.49\linewidth}
    \centering
    \includegraphics[width=\linewidth]{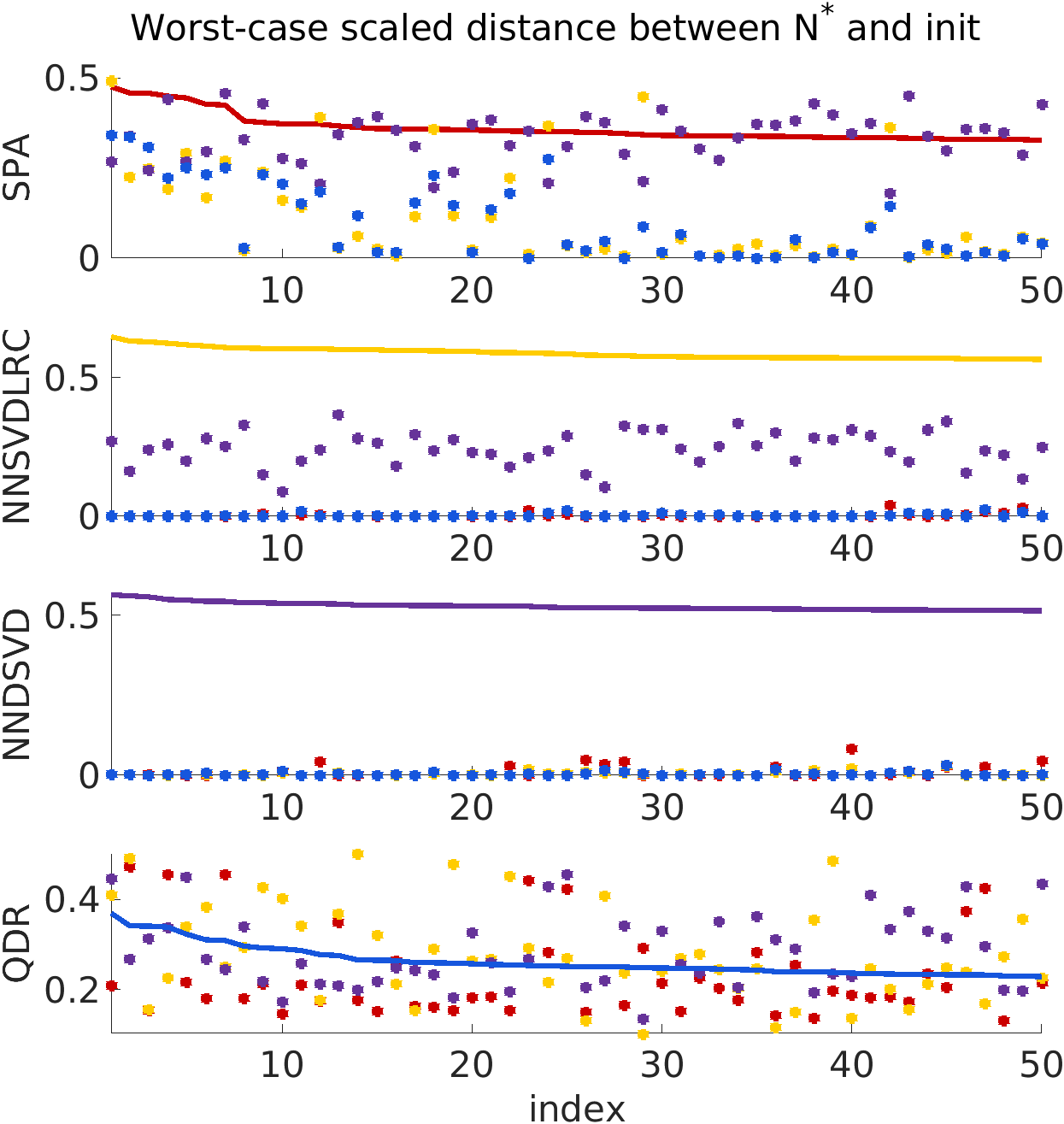}
    \end{subfigure}~
    \begin{subfigure}[t]{0.49\linewidth}
    \centering
    \includegraphics[width=\linewidth]{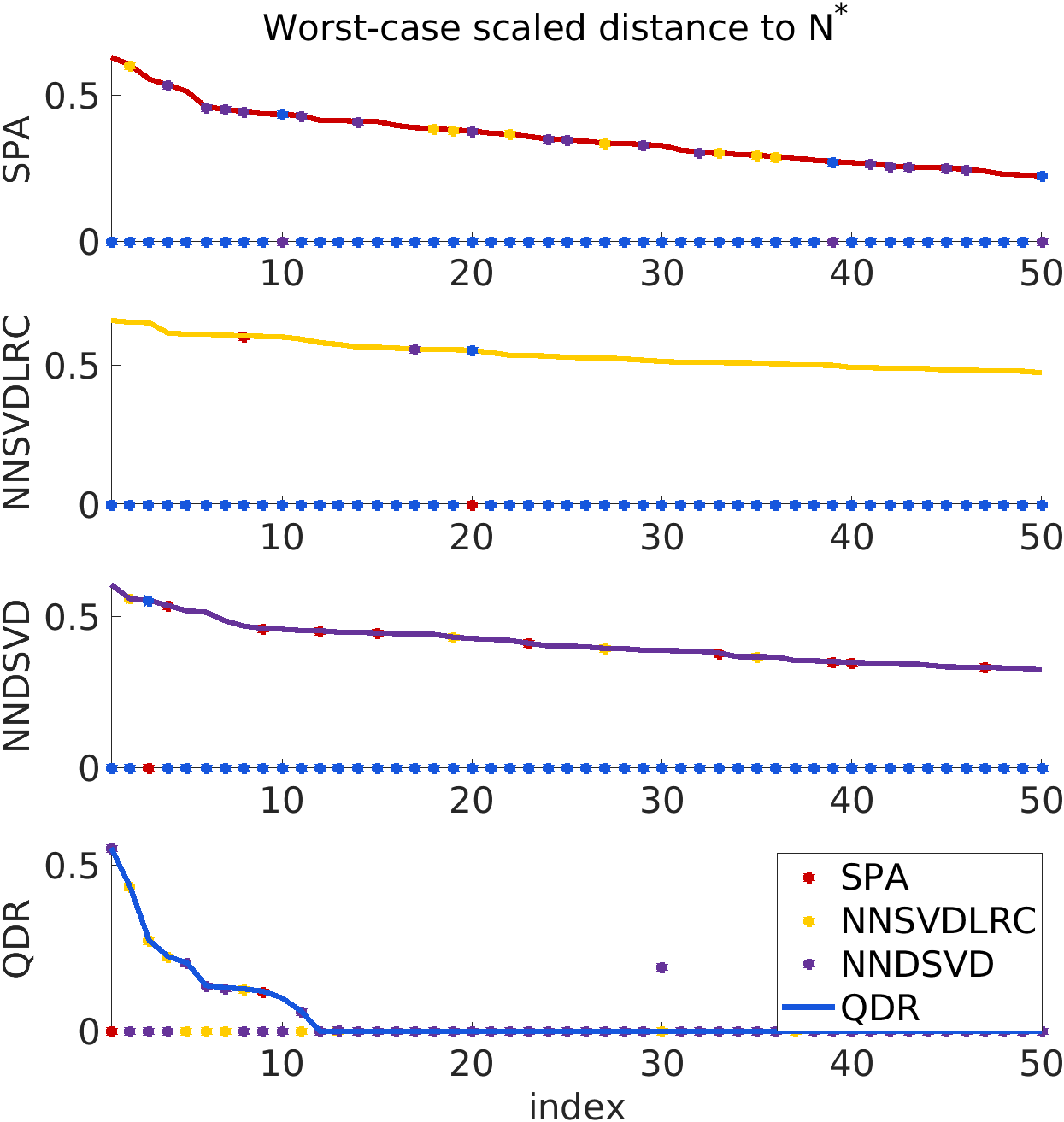}
    \end{subfigure}
    \caption{On the left, the relative distances of the various initialization points to the known optimum $N^*$.
    The 50 worst cases shown for each of the initialization methods.
    On the right, the relative distance in the worst 50 cases after ANLS iterations ($\epsilon = 10^{-5}$, maximum 1000 iterations). }
    \label{fig:comp4x4}
\end{figure}

We also measured how close the starting points are to the true optimum by running Algorithm \ref{alg:quadrant}, as well as its competitors, with no ANLS iterations.
In Figure \ref{fig:comp4x4}, one can see for each initialization method the 50 worst cases, where the output of the algorithm is furthest away from the optimal NMF2 approximation, both for the initialization and the point of convergence of ANLS.

\begin{table}[h]
    \centering
    \begin{tabular}{|c|rrrrr|}
    \hline
 $\log_{10} \delta$  & SPA & NNSVDLRC & NNDSVD & rand & QDR \\ \hline
$\geq -2$ & 79 & 258 & 114 & 1250 & 11 \\
$\in [-5,-2[$ & 10 & 20 & 32 & 15 & 10 \\
$< -5$ &  25268 & 25079 & 25211 & 24092 & 25336 \\
\hline
    \end{tabular}
    \caption{Number of matrices within a distance $\delta:=\| (\cdot) - N^*\|_F/\|N^*\|_F$ from the optimum for the outputs of ANLS with various initialization methods. The ANLS algorithm was ran with tolerance $\varepsilon = 10^{-5}$ and maximum number of iterations 1000.}
    \label{tab:pointofconvergence}
\end{table}

\section{Concluding remarks}\label{sec:conclusions}

We conclude by pointing out some possible directions for future research.

If the original nonnegative matrix $N$ is very large, it might be worth starting from a suboptimal rank-$2$ approximation to save the costs of the SVD. Nonetheless, to use \eqref{eq:DcosD}, one still needs the corresponding vectors $\hat\bu_1$ and $\hat\bv_1$ to be positive and the vectors $\hat\bu_2$ and $\hat\bv_2$ to have both positive and negative elements. 
One possibility is using Lanczos bidiagonalization with a nonnegative starting vector $\hat\bv_1 > 0$ and corresponding vector $\hat\bu_1:=N \hat\bv_1$ that will also be positive as long as $N$ has no zero rows. 
Therefore, if the nonnegative vector $\hat\bv_1$ is a good approximation of the dominant right singular vector of $N$, then the nonnegative vector $\hat\bu_1:=N\hat\bv_1$ will be a good approximation of its left singular vector. Techniques for updating the approximation using ANLS for large sparse matrices, are discussed in \cite{Gillis}.

\section*{Acknowledgements}

We thank Nicolas Gillis for his insightful comments on a preliminary version.

\end{document}